\documentclass[11pt]{amsart}

\usepackage[colorlinks=true, pdfstartview=FitV, linkcolor=blue, 
citecolor=blue]{hyperref}

\usepackage{amssymb,amsmath,amscd}
\usepackage{array}
\usepackage{bbm}
\usepackage{xfrac}
\usepackage{graphicx}
\usepackage{a4wide}
\usepackage{lscape}
\usepackage{tikz}
\usepackage{mathtools}
\tikzset{every picture/.style={line width=0.75pt}} 

\newtheorem{theorem}{Theorem}[section]
\newtheorem{prop}[theorem]{Proposition}
\newtheorem{lemma}[theorem]{Lemma}     
\newtheorem{fact}[theorem]{Fact}
\newtheorem{coro}[theorem]{Corollary}
\theoremstyle{definition}

\newtheorem{remark}[theorem]{Remark}

\newcommand{\ts}{\hspace{0.5pt}}

\newcommand{\RR}{\mathbb{R}\ts}

\newcommand{\ZZ}{\mathbb{Z}}

\newcommand{\NN}{\mathbb{N}}

\newcommand{\QQ}{\mathbb{Q}}

\newcommand{\XX}{\mathbb{X}}

\newcommand{\cA}{\mathcal{A}}

\newcommand{\dd}{\,\mathrm{d}}

\newcommand{\bs}{\boldsymbol}

\newcommand{\Mat}{\mathrm{Mat}}

\renewcommand{\phi}{\varphi}
\renewcommand{\theta}{\vartheta}

\newcommand{\exend}{\hfill$\Diamond$}

\newcommand{\myfrac}[2]{\frac{\raisebox{-2pt}{$#1$}}
      {\raisebox{0.5pt}{$#2$}}}

\numberwithin{equation}{section}

\makeatletter
\renewcommand{\@captionfont}{\small}
\makeatother

\DeclareFontFamily{U}{mathx}{\hyphenchar\font45}
\DeclareFontShape{U}{mathx}{m}{n}{ <5> <6> <7> <8> <9> <10>
   <10.95> <12> <14.4> <17.28> <20.74> <24.88> mathx10 }{}
\DeclareSymbolFont{mathx}{U}{mathx}{m}{n}
\DeclareFontSubstitution{U}{mathx}{m}{n}
\DeclareMathAccent{\widecheck}{0}{mathx}{"71}

\newcommand{\defeq}{\mathrel{\mathop:}=}
\newcommand{\eqdef}{=\mathrel{\mathop:}}

\begin{document}

\title[Correlation functions of  RS sequence]
{Correlation functions of the Rudin--Shapiro sequence}

\author{Jan Maz\'a\v c}
\address{Fakult\"at f\"ur Mathematik, Universit\"at Bielefeld,
  \newline \indent Postfach 100131, 33501 Bielefeld, Germany}
\email{jmazac@math.uni-bielefeld.de}

\begin{abstract} 
In this paper, we show that all odd-point correlation functions of the balanced Rudin--Shapiro sequence vanish and that all even-point correlation functions depend only on a single number, which holds for any weighted correlation function as well. For the four-point correlation functions, we provide a more detailed exposition which reveals some arithmetic structures and symmetries. In particular, we show that one can obtain the autocorrelation coefficients of its topological factor with maximal pure point spectrum among them. 
\end{abstract}

\maketitle

\centerline{Dedicated to the memory of Uwe Grimm}

\bigskip

\section{Introduction}
The \emph{Rudin--Shapiro sequence}, sometimes called Golay--Rudin--Shapiro sequence, is an infinite sequence discovered and studied within the scope of Fourier analysis, independently by Golay, Rudin and Shapiro in 1950 \cite{Golay_1,Golay_2,Rudin,Shapiro}. The original definition of the $n$-th member of the sequence counts the number of two consecutive ones in the binary expansion of $n$. Formally, let $(e_k\ e_{k-1}\ \dots e_2 \ e_1 \ e_0)_2$ be the binary expansion of $n$ and let $b_n$ denote the number of ``$1 \ 1$" in this expansion. We have $b_0=0$ and, for $n\geqslant 1$, the number can  be calculated as 
\[ b_n = \displaystyle \sum_{i=0}^{k-1} e_i e_{i+1}. \]
Then, the $n$-th digit of the (one-sided) Rudin--Shapiro sequence $a_n$ reads
\[a_n \defeq (-1)^{b_n}. \]
This definition also implies that the Rudin--Shapiro sequence is an automatic sequence \cite{Allouche_Shallit},  and provides a possible way to generalise it; see \cite{Q87,AL91} for further details. The standard and generalised Rudin--Shapiro sequences were already studied from different points of view; for the spectral properties of the generalised ones \cite{Chan_Grimm}. 
In particular, one can study the complexity of the Rudin--Shapiro sequence with respect to finite arithmetic progressions contained in the sequence. Konieczny showed that the Rudin--Shapiro sequence is Gowers uniform for any uniformity norm \cite{Konieczny}. 
This result, roughly speaking, demonstrates that the Rudin--Shapiro sequence is not too distant from a random sequence, which is also apparent from its spectral properties \cite{TAO}, as both binary sequences possess absolutely continuous diffraction only. This fact is also visible at the level of the autocorrelation, where, for both sequences, the autocorrelation coefficients vanish except for 0.
In order to understand the statistical behaviour and possible differences better, we exploit the higher-order correlation functions, which are, to some extent, related to the Gowers norms. 
Correlation functions have been widely studied \cite{Domb_Green}, and used in statistical mechanics, and some examples of one-dimensional quasicrystals have been discussed in the corresponding literature, for example \cite{vanE}. It was also shown that one can construct a suitable Ising model such that its Hamiltonian has the Thue--Morse sequence as its ground state once its 4-point correlations are known \cite{GMRE}. 

Recently, Baake and Coons \cite{BC} described the higher-order correlation functions for the Thue--Morse word using renormalisation techniques. We aim to understand the statistical structure of the Rudin--Shapiro word in terms of higher-order correlation functions as well. This may bring a different insight into the type of its long-range order. 

\section{Rudin--Shapiro sequence}
The Rudin--Shapiro (RS) sequence can be obtained via a bi-infinite fixed point of a constant-length substitution over the quaternary alphabet $\cA = \{ 0,1,2,3 \}$, namely 
\[\varrho_{_{\mathrm{RS}}} = \left\{ \begin{array}{rcl}
     0 & \mapsto & 02, \\
     1 & \mapsto & 32, \\
     2 & \mapsto & 01, \\
     3 & \mapsto & 31, \\
\end{array} \right. \]
with any of the legal seeds $2|0$, $2|3$, $1|0$ or $1|3$. The resulting fixed point (under the square of the substitution) defines the \emph{quaternary} RS sequence. Any of them give rise to the \emph{quaternary RS hull} in the standard way. We can further use these quaternary sequences to define the \emph{binary} ones via the mapping  $\phi: \{0,1,2,3\} \longrightarrow \{a,b \}$ defined as 
\[\phi(0) = \phi(2) = a, \qquad \phi(1)=\phi(3)=b.\]
The four bi-infinite fixed points are mapped to four locally indistinguishable binary sequences. Any of them defines, in the usual way, the \emph{binary RS hull}. The binary bi-infinite sequences can be obtained as a fixed point of a \emph{non-local} binary substitution rule.

\begin{lemma}\cite[Lemma 4.9]{TAO}
The four bi-infinite binary Rudin--Shapiro sequences are fixed under the substitution rules
\[ \varrho_{\mathrm{even}} = \left\{\begin{array}{rcl}
     a & \mapsto & aaab, \\
     b & \mapsto & bbba,
\end{array} \right.  \quad \varrho_{\mathrm{odd}} = \left\{\begin{array}{rcl}
     a & \mapsto & aaba, \\
     b & \mapsto & bbab,
\end{array} \right.\]
where $\varrho_{\mathrm{even}}$ and $\varrho_{\mathrm{odd}}$ have to be applied to letters at even and odd positions, respectively. In particular, this rule is non-local, as one needs a reference point to apply the rule. \qed
\end{lemma}

Denote by $\XX_{\mathrm{RS},4}$ the quaternary RS hull and by $\XX_{\mathrm{RS},2}$ the binary one. The mapping $\phi$ induces a continuous mapping commuting with the substitution action, i.e., the following diagram commutes

\begin{equation}
\label{eq:diagram}
\begin{CD}
\XX_{\mathrm{RS},4} @>\varrho_{_{\mathrm{RS}}}>> \XX_{\mathrm{RS},4}\\
@V{\phi}VV @VV{\phi}V\\
\XX_{\mathrm{RS},2} @>\varrho>> \XX_{\mathrm{RS},2}
\end{CD} \quad.
\end{equation}

Moreover, the mapping $\phi$ is invertible. Thus, the hulls are topologically conjugate and even mutually locally derivable (MLD). See \cite[Rem. 4.11]{TAO} for further details. 

From the spectral point of view, the RS sequence is a paradigm of a substitution sequence with an absolutely continuous spectrum, both in the diffraction and in the dynamical sense. For a more detailed discussion on the spectral properties, we refer to \cite{Frank}.

To get deeper inside the structure of the quaternary RS sequence, one can study its two-letter legal subwords \cite{TAO}. In particular, one can define a sliding block map $\chi$,
\begin{align*}
    \chi(01)=\chi(32)=A, & \qquad \chi(02)=\chi(31)=B, \\
    \chi(10) = \chi(23) = C, & \qquad \chi(20)=\chi(13)=D.
\end{align*}

\noindent The quaternary RS substitution then induces a substitution $\varrho^{}_{2}$ on the alphabet $\{A, B, C, D\}$ which reads
\begin{equation}
    \label{eq:subst_induced}
\varrho^{}_{2} = \left\{ \begin{array}{rcl}
     A & \mapsto & BC, \\
     B & \mapsto & BD, \\
     C & \mapsto & AD, \\
     D & \mapsto & AC. \\
\end{array} \right. 
\end{equation}
This substitution rule makes the following diagram with the natural $\ZZ$-action commutative, 

\begin{equation}
\label{eq:diagram_2}
\begin{CD}
\XX_{\mathrm{RS},4} @>S_{\mathrm{RS,4}}>> \XX_{\mathrm{RS},4}\\
@V{\chi}VV @VV{\chi}V\\
\chi(\XX_{\mathrm{RS},4}) @>S^{}_{2}>> \chi(\XX_{\mathrm{RS},4})
\end{CD} \quad.
\end{equation}

The square of the substitution $\varrho^{}_{2}$ possesses two bi-infinite fixed points with starting seeds $C|D$ and $D|B$, respectively. They coincide at all positions but $-1$. As in the RS case, we can use the same trick to obtain a binary sequence, namely the coding $\Tilde{\phi}: \{A,B,C,D\} \longrightarrow \{a,b \}$,
\[\Tilde{\phi}(A) = \Tilde{\phi}(C) = a, \qquad \Tilde{\phi}(B)=\Tilde{\phi}(D)=b.\]
Moreover, the description of this binary sequence as a fixed point of certain, paper--folding-like substitution is possible.\footnote{There is also an interesting connection to space-filling curves arising from this substitution. See Example 5.1.5 in \cite{Allouche_Shallit} for further details.  }
\begin{lemma}
\label{lem:induced}
The bi-infinite binary versions of the fixed points of the substitution $\varrho^{}_{2}$ are fixed under the substitution rules
\begin{equation}
    \label{eq:RS-derived}
 \varrho_{\, \mathrm{2,even}} = \left\{\begin{array}{rcl}
     a & \mapsto & bbab, \\
     b & \mapsto & bbaa,
\end{array} \right.  \quad \varrho_{\, \mathrm{2,odd}} = \left\{\begin{array}{rcl}
     a & \mapsto & baaa, \\
     b & \mapsto & baab,
\end{array} \right.
\end{equation}
where $\varrho_{\mathrm{2,even}}$ and $\varrho_{\mathrm{2,odd}}$ have to be applied to letters at even and odd positions, respectively.

\begin{proof}
Both bi-infinite fixed points $\bs{w}$ of the quaternary substitution $\varrho^{}_{2}$ satisfy $w^{}_{2i} \in \{A,B\}$ and $w^{}_{2i+1} \in \{C,D\}$ for all $i\in\ZZ$. Taking the square of the \eqref{eq:subst_induced}, one gets at even and odd positions
\[ \begin{array}{rlrl}
     A & \mapsto  BDAD, & C & \mapsto  BCAC, \\
     B & \mapsto  BDAC, & D & \mapsto  BCAD. \\
\end{array}  \]
Recalling the coding $\Tilde{\phi}$ gives
the relations. 
\end{proof}
\end{lemma}

\begin{remark}
\label{eq:derived-description}
 Assign weights to each point, namely let $a=-b=1$ and let $m\in\ZZ$. Then, 
\begin{equation*}
w_{4m} = -1, \qquad w_{4m+1} = (-1)^{m+1}, \qquad w_{4m+2} = 1, \qquad w_{4m+3} = (-1)^{m+1}w_m 
\end{equation*}
provide a direct description of the fixed point as a consequence of Lemma \ref{lem:induced}.
\exend
\end{remark}

\begin{remark}
This system possesses a pure point spectrum, as seen by applying Wiener's criterion. Moreover, the comparison with the pure point part of the \emph{dynamical} spectrum of the RS sequence reveals that the dynamical system $(\tilde{\varphi}\chi(\XX_{\mathrm{RS},4}),\ZZ)$ is a topological factor of the RS hull with maximal pure point spectrum \cite{BC, Q_book}. \exend
\end{remark}

\section{Correlation functions}
Due to the MLD relations, we now concentrate on the binary version. Recall that, if one considers the binary RS bi-infinite sequence as a word $\bs{w}$ over the alphabet $\{1,-1\}$, the following recursive relation holds for its letters
\begin{equation}
\label{eq:renorm_letters}
w_{4m+\ell} =  \left\{\begin{array}{rcl}
     w_m, & \mathrm{if} & \ell \in \{0,1\}, \\
     (-1)^{m+\ell}\,w_m, & \mathrm{if} & \ell \in \{2,3\}.
\end{array} \right. 
\end{equation}

We define the $n$-point correlation function $\eta^{(n)}$ as
\begin{equation}
    \label{eq:def_corr_n}
\eta^{(n)}(m_1,m_2,\dots,m_{n-1}) =\lim_{N\to \infty}\myfrac{1}{2N+1} \sum_{i=-N}^{N}w_i\, w_{i+m_1}\, w_{i+m_2}\dots w_{i+m_{n-1}}.
\end{equation}
Note that one can also work with one-sided averages. 
The limit exists due to the (unique) ergodicity of the subshift $\XX_{\mathrm{RS},2} $ via Birkhoff's ergodic theorem. Indeed, one has 
\begin{equation}
\label{eq:def_eta_integral}
\lim_{N\to \infty}\myfrac{1}{2N+1} \sum_{i=-N}^{N}w_i\, w_{i+m_1}\, w_{i+m_2}\dots w_{i+m_{n-1}} \ = \ \int_{\XX_{\mathrm{RS},2}} x_{_0}\, x_{m_1}\, x_{m_2}\dots x_{m_{n-1}} \dd \mu(x) , 
\end{equation}
where  $\mu$ denotes the unique shift-invariant probability measure on the RS shift space, which is the \emph{patch frequency measure} of the subshift. 

In the same way, we define the signed $n$-point correlation function $\theta^{(n)}$, which is well defined for the same reasons as above, by
\begin{equation}
\label{eq:theta-def}
\theta^{(n)}(m_1,m_2,\dots,m_{n-1}) =\lim_{N\to \infty}\myfrac{1}{2N+1} \sum_{i=-N}^{N}(-1)^{i} w_i\, w_{i+m_1}\, w_{i+m_2}\dots w_{i+m_{n-1}}.
\end{equation}

Both $\eta^{(n)}$ and $\theta^{(n)}$ possess several symmetries which play a crucial role in our further investigation. Due to the commutativity of the standard multiplication,  $\eta^{(n)}$ and $\theta^{(n)}$ are invariant under permutations of their arguments. For every permutation $\sigma \in S_{n-1}$, one has
\[\eta^{(n)}(m_1,m_2,\dots,m_{n-1}) = \eta^{(n)}(m_{\sigma(1)},m_{\sigma(2)},\dots,m_{\sigma(n-1)}), \]
\[\theta^{(n)}(m_1,m_2,\dots,m_{n-1}) = \theta^{(n)}(m_{\sigma(1)},m_{\sigma(2)},\dots,m_{\sigma(n-1)}). \]
Further,  in \eqref{eq:def_eta_integral}, one can shift the summation index and omit several terms. Indeed, for example, for $\eta^{(n)}$ one has for any fixed $t\in \ZZ$
\begin{align*}
\eta^{(n)}(m_1,m_2,\dots,m_{n-1}) {}&=\lim_{N\to \infty}\myfrac{1}{2N+1} \sum_{i=-N}^{N}w_i\, w_{i+m_1}\, w_{i+m_2}\dots w_{i+m_{n-1}} \\
{}&=\lim_{N\to \infty}\myfrac{1}{2N+1} \sum_{i=-N-t}^{N-t}w_{i+t}\, w_{i+t+m_1}\, \dots w_{i+t+m_{n-1}}\\
 {}&= \lim_{N\to \infty}\myfrac{1}{2N+1}\left( \sum_{i=-N-t}^{-N+t}w_{i+t}\, w_{i+t+m_1}\, \dots w_{i+t+m_{n-1}} \right. \\
 {}& \qquad \qquad \qquad \qquad \left. + \sum_{i=-N+t}^{N-t}w_{i+t}\, w_{i+t+m_1}\, \dots w_{i+t+m_{n-1}}\right) \\
 {}&= \ 0 + \lim_{N\to\infty}\myfrac{1}{2N-2|t|+1}\sum_{i=-N+t}^{N-t}w_{i+t}\, w_{i+t+m_1}\, \dots w_{i+t+m_{n-1}}. 
\end{align*}
In the calculation, we used that $\left|\sum_{i=-N-t}^{-N+t}w_{i+t}\, w_{i+t+m_1}\, \dots w_{i+t+m_{n-1}} \right|\leqslant 2|t|+1$ for all $N$ and $\myfrac{2N-2|t|+1}{2N+1} \xrightarrow{N \to \infty} 1$. This property of the correlation function becomes nicer if one tracks all relative positions of the elements of the sums, i.e., if one includes the 0 term. That is, we write
\[\langle 0, m_1,m_2,\dots,m_{n-1} \rangle \defeq \eta^{(n)}(m_1,m_2,\dots,m_{n-1}) . \]
Now, the symmetry described above is nothing but a translation symmetry within $\ZZ^n$ in the direction $(1,1,\dots, 1)$, so
\[ \langle 0, m_1,m_2,\dots,m_{n-1} \rangle = \langle \ell, m_1+\ell,m_2+\ell,\dots,m_{n-1}+\ell \rangle \quad \mbox{for every } \ell \in \ZZ. \]
To rewrite the last expression again as the correlation function, one needs to obtain at least one zero in the expression $\langle \ell, m_1+\ell,m_2+\ell,\dots,m_{n-1}+\ell \rangle $. This corresponds to the choices $\ell \in \{0,-m_1, -m_2, \dots, - m_{n-1} \}$. Then, we obtain a ``shadow" of this $\ZZ^n$-translation in the $\ZZ^{n-1}$-dimensional coordinate space, namely 
\[ \eta^{(n)}(m_1,m_2,\dots,m_{n-1}) = \eta^{(n)}(-m_1,m_2-m_1,\dots,m_{n-1}-m_1). \]

Together with the permutation-invariance of the correlation functions, we can attach to each point its orbit under these two symmetries. It turns out that there are at most $n \cdot (n-1)!$ elements in each orbit. 

Moreover, these two symmetries permit a restriction to positive integers, namely to $\ZZ^{n-1}_{\geqslant 0}$. We show that every element of $\ZZ^{n-1}\backslash \ZZ^{n-1}_{\geqslant 0}$ lies in an orbit of an element from $\ZZ^{n-1}_{\geqslant 0}$. Let $(m_1,m_2,\dots,m_{n-1}) \in \ZZ^{n-1}\backslash \ZZ^{n-1}_{\geqslant 0}$ and let us assume that $|m_1|\geqslant |m_2| \geqslant \dots \geqslant |m_{n-1}|$ and let $m_1<0$. This can always be achieved via a suitable permutation. Then, since 
\[ \eta^{(n)}(m_1,m_2,\dots,m_{n-1}) = \eta^{(n)}(-m_1,m_2-m_1,\dots,m_{n-1}-m_1), \]
$(-m_1,m_2-m_1,\dots,m_{n-1}-m_1)$ is a positive vector that can be used for the calculations. 

Since we consider a word over the binary alphabet $\{-1,1\}$ we also get the so-called \emph{cancellation property}. Because $w_i^2=1$ for every $i$, we have for any $n\geqslant2$
\[ \eta^{(n+2)}(m_1,m_2,\dots,M,\dots,M,\dots,m_{n-1}) =\eta^{(n)}(m_1,m_2,\dots,m_{n-1}), \]
so a reduction is possible whenever two entries agree. 
Let us summarise as follows. 
\begin{prop}
\label{prop:symmetries}
Let $\bs{w}$ be a bi-infinite word over the alphabet $\mathcal{A}$ and let $\eta^{(n)}$ and $\theta^{(n)}$ denote the $n$-point correlations functions of the word $\bs{w}$ according to \eqref{eq:def_eta_integral} and \eqref{eq:theta-def}. Further, let us suppose that $\eta^{(n)}$ and $\theta^{(n)}$ are well defined. Then, they possess the following symmetries. 
\begin{enumerate}
    \item Invariance under the action of $S_{n-1}$, namely for every $\sigma \in S_{n-1}$
    \[\eta^{(n)}(m_1,m_2,\dots,m_{n-1}) = \eta^{(n)}(m_{\sigma(1)},m_{\sigma(2)},\dots,m_{\sigma(n-1)}), \]
\[\theta^{(n)}(m_1,m_2,\dots,m_{n-1}) = \theta^{(n)}(m_{\sigma(1)},m_{\sigma(2)},\dots,m_{\sigma(n-1)}), \]
    \item Invariance under a higher-dimensional translation, namely
    \[ \eta^{(n)}(m_1,m_2,\dots,m_{n-1}) = \eta^{(n)}(-m_1,m_2-m_1,\dots,m_{n-1}-m_1). \]
    \[ \theta^{(n)}(m_1,m_2,\dots,m_{n-1}) = (-1)^{m_1}\,\theta^{(n)}(-m_1,m_2-m_1,\dots,m_{n-1}-m_1). \]
\end{enumerate}
The orbit of any point in $\ZZ^{n-1}$ under these two operations has at most $n!$ elements. 

If the word $\bs{w}$ is over the binary alphabet $\{-1,1\}$, the correlation functions with $n\geqslant2$ have the \emph{cancellation property}, i.e. 
\[ \eta^{(n+2)}(m_1,m_2,\dots,M,\dots,M,\dots,m_{n-1}) =\eta^{(n)}(m_1,m_2,\dots,m_{n-1}), \]
\[ \pushQED{\qed} 
\theta^{(n+2)}(m_1,m_2,\dots,M,\dots,M,\dots,m_{n-1}) =\theta^{(n)}(m_1,m_2,\dots,m_{n-1}). \qedhere \popQED \]

\end{prop}

\begin{remark}
For the autocorrelation function $\eta^{(2)}$, the second property is nothing but a~mirror symmetry with respect to the origin. In other words,
\[ \eta^{(2)}(m) {} = {} \eta^{(2)}(-m) \]
holds for every $m \in \ZZ$.
This mirror symmetry is generally not satisfied for higher-order correlation functions. \exend
\end{remark}

\section{Renormalisation relations}
This section aims to present the renormalisation structure of the correlation functions of the RS sequence and investigate their behaviour, based on the recursive relation \eqref{eq:renorm_letters} and the ergodicity of the subshift. The renormalisation equations can be derived by splitting the summation into four terms; in each, the summation index goes over the residues modulo 4. We demonstrate this procedure for one particular case, and its generalisation will be obvious. 

Let us start with the three-point correlation function $\eta^{(3)}$. We wish to determine its value at $(4m_1+2,4m_2+3)$ using the values of $\eta^{(3)}$ and $\theta^{(3)}$ at  $(m_1+r_1,m_2+r_2)$ with $r_i\in\{0,1,2,3\}$.

\begin{align*}
    \eta^{(3)}(4m_1 +&{}2,4m_2 +3) = \lim_{N\to \infty}\myfrac{1}{2N+1} \sum_{i=-N}^{N}w_i\, w_{i+4m_1+2}\, w_{i+4m_2+3} \\
    {}&=\lim_{N\to \infty}\myfrac{1}{8N\! + \!1} \sum_{i=-4N}^{4N}w_i\, w_{i+4m_1+2}\, w_{i+4m_2+3} \\
    {}&= \begin{multlined}[t] 
            \lim_{N\to \infty}\myfrac{1}{8N\! + \!1} \sum_{j=-N}^{N}\left(w_{4j}\, w_{4j+4m_1+2}\, w_{4j+4m_2+3} +w_{4j+1}\, w_{4j+4m_1+3}\, w_{4j+4m_2+4}\right. \\
            \left. + w_{4j+2}\, w_{4j+4m_1+4}\, w_{4j+4m_2+5} +w_{4j+3}\, w_{4j+4m_1+5}\, w_{4j+4m_2+6} \right)
        \end{multlined}\\
    {}&=\begin{multlined}[t]
            \lim_{N\to \infty}\myfrac{1}{8N\! + \!1} \sum_{j=-N}^{N}\!\bigl(\!-\!(-1)^{m_1+m_2}w_{j} w_{j+m_1}w_{j+m_2}\! -\!(-1)^{m_1+j}w_{j} w_{j+m_1} w_{j+m_2+1} \\
             + (-1)^{j}w_{j}\, w_{j+m_1+1}\, w_{j+m_2+1} + (-1)^{m_2}w_{j}\, w_{j+m_1+1}\, w_{j+m_2+1} \bigr)
        \end{multlined}\\
    {}&=\begin{multlined}[t]
            \myfrac{1}{4}\left(-(-1)^{m_1+m_2}\eta^{(3)}(m_1,m_2)-(-1)^{m_1} \theta^{(3)}(m_1,m_2\!+\!1)+\theta^{(3)}(m_1\!+\!1,m_2\!+\!1) \right.\\
            \left. +(-1)^{m_2}\eta^{(3)}(m_1\!+\!1,m_2\!+\!1)\right).
        \end{multlined}
\end{align*}
This calculation reveals a general (and useful!) fact which holds for an arbitrary correlation function of the RS sequence. Namely that, independently on the left hand side of the renormalisation equation, on the right hand side one gets arguments containing $m_i$ or $m_i+1$ only (and does not contain $m_i+2$ and $m_i+3$). This observation is a simple consequence of the recursion \eqref{eq:renorm_letters}, but it has some far-reaching consequences. 

The renormalisation equations form an infinite set of linear equations in infinitely many variables which can be split into two parts --- the \emph{self-consistent} part and the \emph{recursive} part. The first one consists of equations that cannot be simplified via the renormalisation equations. This self-consistent part closes on itself. The above helps us determine this \emph{finite} self-consistent part.
The recursive part of the renormalisation equations is then entirely determined by the solution of the self-consistent part. This shows that the dimension of the solution space of the renormalisation equations is \emph{finite}. 

Let us start the self-consistent part of the renormalisation equations for $n=2$ and $n=3$. It turns out that these two solutions are sufficient for determining the correlation functions for arbitrary $n$. Setting $m_i=0$ for all $i$ on the RHS of the renormalisation equations leads to the self-consistent part with arguments $(0),(1)$, and $(0,0), (0,1), (1,1)$, respectively.

For the self-consistent part of the 2-point correlation function, one obtains
\begin{align*}
\label{eq:s_c_2point}
\eta^{(2)}(0) &=  \eta^{(2)}(0), & \theta^{(2)}(0) &=  0, \\
\eta^{(2)}(1) &=  \myfrac{1}{4}\left( \theta^{(2)}(0) - \theta^{(2)}(1)\right), & \theta^{(2)}(1) &=  \myfrac{1}{4}\left( -\theta^{(2)}(0) + \theta^{(2)}(1)\right). 
\end{align*}
The dimension of the solution space is one, and the unique solution is fully determined from $\eta^{(2)}(0)$. This value can be calculated directly from the definition, which gives $\eta^{(2)}(0) =1$. These results are well known and can be found together with the set of all renormalisation equations for the autocorrelation function in \cite[Sec. 10.2.]{TAO}. We state the result as follows. 

\begin{lemma}
The autocorrelation coefficients $\eta^{(2)}$ of the signed Dirac comb of the binary Rudin--Shapiro sequence exist for all $m\in\ZZ$ and are given by  $\eta^{(2)}(m) = \delta_{m,0}$. \qed
\end{lemma}

This result excludes the RS sequence to be pure point diffractive due to Wiener's criterion \cite[Prop. 8.9.]{TAO}. Employing the explicit relation for $\eta^{(2)}(m)$, we can derive the autocorrelation measure $\gamma = \sum_{m\in \ZZ} \eta^{(2)}\!(m)\,\delta_m$, and its Fourier transform $\widehat{\gamma}$ being the Lebesgue measure. This shows that the balanced RS word possesses an absolutely continuous spectrum only. It is perhaps unexpected that there is no difference between the RS sequence (a fixed point of a primitive substitution) and a binary Bernoulli sequence in $\{\pm1\}$ with equal probabilities on the level of autocorrelation coefficients. Therefore, it is worth studying the higher-order correlations to understand the differences between these two structures. For the Bernoulli sequence, all correlations vanish. We expect this to be different for the RS sequence. 

Let us move to the case of 3-point correlations. The complete set of renormalisation equations can be found in the Appendix. The self-consistent part of this system simply reads
\begin{align*}
\eta^{(3)}(0,0) {}&= \myfrac{1}{2}\, \eta^{(3)}(0,0) , \\
\eta^{(3)}(0,1) {}&= \myfrac{1}{4} \left(\eta^{(3)}(0,0) + \eta^{(3)}(0,1)\right),  \\ 
\eta^{(3)}(1,1) {}&= \myfrac{1}{4} \left(2\eta^{(3)}(0,0) + \theta^{(3)}(0,0) - \theta^{(3)}(1,1) \right),  \\ 
\theta^{(3)}(0,0) {}&= 0, \\
\theta^{(3)}(0,1) {}&= \myfrac{1}{4} \left(\eta^{(3)}(0,0)- 2\theta^{(3)}(0,0) - \eta^{(3)}(0,1) \right), \\
\theta^{(3)}(1,1) {}&= \myfrac{1}{4} \left(\theta^{(3)}(0,0) + \theta^{(3)}(1,1)\right),
\end{align*}
and has the trivial solution $\eta^{(3)}(k,\ell) = \theta^{(3)}(k,\ell) = 0$ only, as one can easily see. Thus, we can profit from these two results and  the properties of the correlation functions as follows.  

\begin{prop}
The $n$-point correlation functions $\eta^{(n)}, \theta^{(n)}$ of the signed Dirac comb of the binary Rudin--Shapiro sequence exist for all $n\in\ZZ$. They are all determined by the value of $\eta^{(n)}(0,\dots,0)$. For odd $n$, the correlation functions $\eta^{(n)}, \theta^{(n)}$ vanish. 
\begin{proof}
The existence follows from the unique ergodicity of the system and Birkhoff's ergodic theorem. The calculation above shows that the functions in self-consistent part of the renormalisation equations for 3-point correlations vanishes. The recursive structure of the problem implies that all three-point correlation functions are zero. 
Consider now any $(2k+1)$-point correlation function and any point $x \in \ZZ^{2k}$ which belongs to the self-consistent part of the renormalisation equations. As discussed above, we have $x\in \{0,1\}^{2k}$. From the pigeonhole principle and the cancellation property, it follows that the self-consistent part of any correlation function can always be reduced to the 2- or 3-point one, which gives the one-dimensional solution and, moreover, implies that the odd-$n$ correlation function vanishes. 
\end{proof}
\end{prop}
As a direct consequence of the previous proof, we obtain the values of the correlation functions $\eta^{(n)}, \ \theta^{(n)}$ on the vertices of the $(n{-}1)$-dimensional hypercube. 

\begin{coro}
Let $(r_1, \ r_2, \dots, r_{n-1}) \in \{0,1\}^{n-1}$ and denote $r = r_1 +r_2 + \cdots + r_{n-1}$. Then, for the correlation functions $\eta^{(n)},\ \theta^{(n)}$ of the binary Rudin--Shapiro word, one has  
\begin{equation}
    \eta^{(n)}(r_1, \ r_2, \dots, r_{n-1}) =  \left\{\begin{array}{rl}
     0, & \mbox{ if $n$ or $r$ are odd}, \\[5pt]
     1, & \mbox{ if $n$  and $r$ are even},
\end{array} \right. 
\end{equation}
and $\theta^{(n)}(r_1, \ r_2, \dots, r_{n-1}) = 0$ for all $n$ and for all \/ $r_i$.
\qed
\end{coro}

We can harvest the underlying idea in the proof even further. Namely, we can generalise the formula above to arbitrary binary substitution with weights $\pm 1$. 

\begin{prop}
Let $\bs{w}$ be an (bi)-infinite word over binary alphabet $\{\pm1\}$, and assume that the letter frequencies exist. Let the correlation functions $\eta^{(n)}$ exist for arbitrary $n$. Let further  $(r_1, \ r_2, \dots, r_{n-1}) \in \{0,1\}^{n-1}$ and set  $r = r_1 +r_2 + \cdots + r_{n-1}$. Then, the correlation function $\eta^{(n)}$, $n\geqslant3$, can only reach three values on these vertices, namely

\begin{equation}
    \eta^{(n)}(r_1, \ r_2, \dots, r_{n-1}) =  \left\{\begin{array}{rl}
     \nu_{1}-\nu_{-1}, & \mbox{ if $n$ is odd}, \\[5pt]
     \eta^{(2)}(0), & \mbox{ if $n$ and $r$ are even},\\[5pt]
     \eta^{(2)}(1), & \mbox{ if $n$ is even and $r$ is odd},
\end{array} \right.
\end{equation}
where $\nu_{a}$ denotes the letter frequency of letter $a$ in the word $\bs{w}$. 
\begin{proof}
Let us start with $n$ even. The pigeonhole principle implies that one can profit from the cancellation property. Namely, the summands in the defining sum \eqref{eq:def_corr_n} can be reduced to one of the following products --- in case of $r$ being even $w_iw_i$ or $w_{i+1}w_{i+1}$, and $w_i w_{i+1}$ otherwise. The resulting limits then give the $\eta^{(2)}(0)$ in the first two cases and $ \eta^{(2)}(1)$ in the last one. 
If $n$ is odd, the products can be reduced to one of the following four triples: $w_iw_iw_i$,  $w_{i+1}w_{i+1}w_{i+1}$, $w_{i+1}w_{i+1}w_{i}$ and $w_iw_i w_{i+1}$. Nevertheless, we still can once again use that fact that $w_i^2=1$ and reduce the products to either $w_i$ or $w_{i+1}$. Then, we get $\eta^{(n)}(r_1, \ r_2, \dots, r_{n-1}) = \lim_{N\to \infty}\frac{1}{2N+1}\sum_{i=-N}^{N}w_i= \lim_{N\to \infty}\frac{1}{2N+1}\sum_{i=-N}^{N}w_{i+1} = \nu_{1}-\nu_{-1}$, which completes the proof. 
\end{proof}
\end{prop}

We have already determined the values of arbitrary correlation function $\eta^{(n)}$ for the case of balanced weights, i.e., $w_i\in\{\pm 1\}$ for all $i\in \ZZ$. Now, we can extend the result to the case of $n$-point correlation function for the RS word with general weights given by an arbitrary weight function $f : \{\pm1 \} \to \RR$. The  $n$-point correlation function with weight $f$ is defined as

\[ \eta_f^{(n)}(m_1,\dots, m_{n-1}) \defeq \int_{\XX_{\mathrm{RS},2}} f(x_0)f(x_{m_1}) \dots f(x_{m_{n-1}}) \dd \mu(x). \]

\noindent Further, we can follow the steps in \cite[Sec. 5]{BC} and define two values playing a key role in the description of $\eta_f^{(n)}$, namely the expectation 
\[ \mathbb{E}(f) = \int_{\XX_{\mathrm{RS},2}} f(x_0)\dd \mu(x) = \myfrac{f(1)+f(-1)}{2}, \quad \mbox{and}\quad  h_f\eqdef \myfrac{f(1)-f(-1)}{2}. \]
Rephrasing all arguments given by Baake and Coons, which are true not only in the case of Thue--Morse word, but in general binary ones with $\nu_{1}-\nu_{-1}=0$  as well (see Proposition~5.1. and the following discussion in \cite{BC}), one obtains the desired result. 

\begin{theorem}
For any $n\geqslant2$ and for any $f : \{\pm1 \} \to \RR$, the $n$-point correlation function of $f$-weighted binary Rudin--Shapiro word can be calculated from the balanced correlations $\eta^{(n)}$. In particular, the functions $\eta_f^{(n)}$ are determined by the single value $\eta^{(n)}(0,\dots,0)$.
\qed
\end{theorem}

\section{Matrix representation}
In the previous section, we introduced the renormalisation equations and (in the Appendix) gave the list for 3- and 4-point correlation functions. It would be helpful to derive a general formula for an arbitrary renormalisation equation for arbitrary $n$ as Baake and Coons did in~\cite{BC} for the Thue--Morse sequence where they also introduced a suitable matrix formalism. In what follows, we would like to do the same procedure for the RS case. For better readability, we omit the upper index in the notation. 

Then, the matrix version of the renormalisation equations reads 
\begin{equation*}
    \begin{pmatrix}
    \eta(4m) \\ \eta(4m+1) \\ \eta(4m+2) \\ \eta(4m+3) \\ \theta(4m) \\ \theta(4m+1) \\ \theta(4m+2) \\ \theta(4m+3)
    \end{pmatrix} = \myfrac{1}{4} \begin{pmatrix}
    2+2(-1)^m & 0 & 0 & 0 & 0 & 0 & 0 & 0 \\
    1-(-1)^m & 0 & 0 & 0 & (-1)^m & -1 & 0 & 0 \\
    0 & 0 & 0 & 0 & 0 & 0 & 0 & 0 \\
    0 & 1+(-1)^m & 0 & 0 & -(-1)^m & 1 & 0 & 0 \\
    0 & 0 & 0 & 0 & 0 & 0 & 0 & 0 \\
    1-(-1)^m & 0 & 0 & 0 & -(-1)^m & 1 & 0 & 0 \\
    0 & 0 & 0 & 0 & 2(-1)^m & 2 & 0 & 0 \\
    0 & 1+(-1)^m & 0 & 0 & -(-1)^m & 1 & 0 & 0 \\
    \end{pmatrix}
    \begin{pmatrix}
    \eta(m) \\ \eta(m+1) \\ \eta(m+2) \\ \eta(m+3) \\ \theta(m) \\ \theta(m+1) \\ \theta(m+2) \\ \theta(m+3)
    \end{pmatrix}.
\end{equation*}
There is a disadvantage in this description consisting in dealing with the \emph{non-locality} of our substitution, which shows here via the terms with $(-1)^m$. Therefore, we must double the dimension and treat the odd and even positions separately. To do so, let us introduce $\bs{\eta}(8m+r) \in \RR^{8+8}$ with entries 
\[
\bs{\eta}(8m+r)_k =   \left\{\begin{array}{rcl}
     \eta(8m+r+i), & \mathrm{if} & k=i, \\[10pt]
     \theta(8m+r+i), & \mathrm{if} & k=8+i.
\end{array} \right.
\]

\noindent Then, we can rewrite the renormalisation equations for the autocorrelation as 
\[\bs{\eta}(8m+r) = \bs{B}_{(r)} \ \bs{\eta}(m).  \]
To proceed further, we introduce a set of 8 integer matrices $M_{(i,j)}$ and $N_{(i,j)}$ of dimension eight, which will later naturally appear in the description of the matrices $\bs{B}$. 

\begin{align*}
  M_{(0,0)} =  \begin{pmatrix} 0 & 0 & 0 & 0 & 0 & 0 & 0 & 0 \\ -1 & 0 & 0 & 0 & 0 & 0 & 0 & 0 \\2 & 0 & 0 & 0 & 0 & 0 & 0 & 0 \\-1 & 0 & 0 & 0 & 0 & 0 & 0 & 0 \\ 0 & 0 & 0 & 0 & 0 & 0 & 0 & 0 \\0 & 1 & 0 & 0 & 0 & 0 & 0 & 0 \\ 0 & -2 & 0 & 0 & 0 & 0 & 0 & 0 \\0 & 1 & 0 & 0 & 0 & 0 & 0 & 0   \end{pmatrix}, & \qquad  M_{(0,1)} = \begin{pmatrix} 0 & 0 & 0 & 0 & 0 & 0 & 0 & 0 \\ 1 & 0 & 0 & 0 & 0 & 0 & 0 & 0 \\0 & 0 & 0 & 0 & 0 & 0 & 0 & 0 \\0 & -1 & 0 & 0 & 0 & 0 & 0 & 0 \\ 0 & 0 & 0 & 0 & 0 & 0 & 0 & 0 \\0 & 1 & 0 & 0 & 0 & 0 & 0 & 0 \\ 0 & 0 & 0 & 0 & 0 & 0 & 0 & 0 \\0 & 0 & -1 & 0 & 0 & 0 & 0 & 0   \end{pmatrix}, 
\end{align*}
\begin{align*}
  M_{(1,0)} = \begin{pmatrix} 0 & 0 & 0 & 0 & 0 & 0 & 0 & 0 \\ 1 & 0 & 0 & 0 & 0 & 0 & 0 & 0 \\0 & 0 & 0 & 0 & 0 & 0 & 0 & 0 \\-1 & 0 & 0 & 0 & 0 & 0 & 0 & 0 \\ 0 & 0 & 0 & 0 & 0 & 0 & 0 & 0 \\0 & -1 & 0 & 0 & 0 & 0 & 0 & 0 \\ 0 & 0 & 0 & 0 & 0 & 0 & 0 & 0 \\0 & 1 & 0 & 0 & 0 & 0 & 0 & 0   \end{pmatrix}, & \qquad  M_{(1,1)} = \begin{pmatrix} 2 & 0 & 0 & 0 & 0 & 0 & 0 & 0 \\ 1 & 0 & 0 & 0 & 0 & 0 & 0 & 0 \\0 & 0 & 0 & 0 & 0 & 0 & 0 & 0 \\0 & 1 & 0 & 0 & 0 & 0 & 0 & 0 \\ 0 & 2 & 0 & 0 & 0 & 0 & 0 & 0 \\0 & 1 & 0 & 0 & 0 & 0 & 0 & 0 \\ 0 & 0 & 0 & 0 & 0 & 0 & 0 & 0 \\0 & 0 & 1 & 0 & 0 & 0 & 0 & 0   \end{pmatrix}, 
\end{align*}
\begin{align*}
  N_{(0,0)} = \begin{pmatrix} 2 & 0 & 0 & 0 & 0 & 0 & 0 & 0 \\ -1 & 0 & 0 & 0 & 0 & 0 & 0 & 0 \\0 & 0 & 0 & 0 & 0 & 0 & 0 & 0 \\0 & 1 & 0 & 0 & 0 & 0 & 0 & 0 \\ 0 & -2 & 0 & 0 & 0 & 0 & 0 & 0 \\0 & 1 & 0 & 0 & 0 & 0 & 0 & 0 \\ 0 & 0 & 0 & 0 & 0 & 0 & 0 & 0 \\0 & 0 & -1 & 0 & 0 & 0 & 0 & 0   \end{pmatrix},
  & \qquad  N_{(0,1)} = \begin{pmatrix} 0 & 0 & 0 & 0 & 0 & 0 & 0 & 0 \\ 0& -1 & 0 & 0 & 0 & 0 & 0 & 0 \\0 & 0 & 0 & 0 & 0 & 0 & 0 & 0 \\0 & 1 & 0 & 0 & 0 & 0 & 0 & 0 \\ 0 & 0 & 0 & 0 & 0 & 0 & 0 & 0 \\0 & 0 & -1 & 0 & 0 & 0 & 0 & 0 \\ 0 & 0 & 0 & 0 & 0 & 0 & 0 & 0 \\0 & 0 & 1 & 0 & 0 & 0 & 0 & 0   \end{pmatrix},
\end{align*}
\begin{align*}
  N_{(1,0)} = \begin{pmatrix} 0 & 0 & 0 & 0 & 0 & 0 & 0 & 0 \\ -1 & 0 & 0 & 0 & 0 & 0 & 0 & 0 \\0 & 0 & 0 & 0 & 0 & 0 & 0 & 0 \\0 &-1 & 0 & 0 & 0 & 0 & 0 & 0 \\ 0 & 0 & 0 & 0 & 0 & 0 & 0 & 0 \\0 & 1 & 0 & 0 & 0 & 0 & 0 & 0 \\ 0 & 0 & 0 & 0 & 0 & 0 & 0 & 0 \\0 & 0 & 1 & 0 & 0 & 0 & 0 & 0  \end{pmatrix}, & \qquad  N_{(1,1)} = \begin{pmatrix} 0 & 0 & 0 & 0 & 0 & 0 & 0 & 0 \\ 0& 1 & 0 & 0 & 0 & 0 & 0 & 0 \\0 & 2 & 0 & 0 & 0 & 0 & 0 & 0 \\0 & 1 & 0 & 0 & 0 & 0 & 0 & 0 \\ 0 & 0 & 0 & 0 & 0 & 0 & 0 & 0 \\0 & 0 & 1 & 0 & 0 & 0 & 0 & 0 \\ 0 & 0 & 2 & 0 & 0 & 0 & 0 & 0 \\0 & 0 & 1 & 0 & 0 & 0 & 0 & 0   \end{pmatrix}. \\ 
\end{align*}
Using these matrices, the matrix $\bs{B}_{(0)} \in \Mat(16,\QQ)$ decomposes into a block matrix,
\begin{equation}
    4\bs{B}_{(0)} =\begin{pmatrix}N_{(0,0)} & N_{(0,1)} \\[10pt] N_{(1,0)} & N_{(1,1)} \end{pmatrix} + \begin{pmatrix}M_{\tau(0,0)} & M_{\tau(0,1)} \\[10pt] M_{\tau(1,0)} & M_{\tau(1,1)} \end{pmatrix}, 
\end{equation}
where $\tau$ is defined for every $(i,j)\in\{0,1\}^2$ as $\tau(i,j) = (i+1,j+1) \mod 2$. All remaining matrices $\bs{B}_{(r)}$ can be obtained from $\bs{B}_{(0)}$ via a power of some transformation $\mathcal{J}$. This mapping naturally encodes the rearrangement of the renormalisation equations done by going from $\bs{\eta}(8m)$ to $\bs{\eta}(8m+1)$.  It shifts the $i$-th row of the matrix to the position $i-1$, and the first row becomes the last with a column shift.

In terms of matrices, we can rewrite it as follows
\begin{equation}
\label{eq:transf_J}
\mathcal{J}(X)  = R_1\cdot X + L_7\cdot X\cdot S \quad \mbox{for any $X\in\mathrm{span}\{M_{(i,j)},N_{(i,j)}\}$}, 
\end{equation}
with $S$ being the permutation matrix (in the column notation, meaning that $s_{ij} = 1$ iff $j=\pi(i)$) of the permutation $\pi=(1\,3\,5\,7)(2\,4\,6\,8)$ and for $m,n \leqslant 7$

\[\bigl(R_m\bigr)_{i,\,j} = \delta_{i,\, j-k}, \qquad \bigl(L_n\bigr)_{i,\,j} = \delta_{i-n,\,j}.  \]

Note that $S^4 = \mathrm{Id}$, and that the matrices $R_1$ and $L_7$ are of rank $7$ and $1$ respectively, implying  $R_1^8 =O$ and $L_7^2=O$ with $O$ standing for the zero matrix. These properties enable one to explicitly write down the powers of $\mathcal{J}$. 

\begin{remark}
For the transformation $\mathcal{J}$ defined in \eqref{eq:transf_J} it is not hard to show that 
\[\mathcal{J}^k(X) = R_{k}\cdot X + L_{8-k}\cdot X \cdot S \]
holds for every $k\in\{1,2,\dots 7\}$ and for every $X\in\mathrm{span}\{M_{(i,j)},N_{(i,j)}\}$ 
with matrices $R_m$ and $L_n$ defined above. 
\exend
\end{remark}

With this transformation, we can  get the desired matrices $\bs{B}_{(r)}$ for $r\in\{1,2,\dots 7\}$ as 
\begin{equation}
    4\bs{B}_{(r)} = 4\mathcal{J}^r\bs{B}_{(0)} = \begin{pmatrix}\mathcal{J}^r(N_{(0,0)}) & \mathcal{J}^r(N_{(0,1)}) \\[10pt] \mathcal{J}^r(N_{(1,0)}) & \mathcal{J}^r(N_{(1,1)}) \end{pmatrix} + \begin{pmatrix}\mathcal{J}^r(M_{\tau(0,0)}) & \mathcal{J}^r(M_{\tau(0,1)}) \\[10pt] \mathcal{J}^r(M_{\tau(1,0)}) & \mathcal{J}^r(M_{\tau(1,1)}) \end{pmatrix}.
\end{equation}

In order to describe the matrices for the higher-order correlation functions, we have to generalise the notion of vectors $\bs{\eta}$. In particular, for 3-point correlations, we employ the vectors $\bs{\eta}(8m_1+r_1, 8m_2+r_2) \in\RR^{8^2+8^2}$ with entries
\begin{equation} 
\label{eq:eta_def}
\bs{\eta}(8m_1+r_1,\ 8m_2+r_2)_k =   \left\{\begin{array}{rcl}
     \eta(8m_1+r_1+i,\ 8m_2+r_2+j), & \mathrm{if} & k=8i+j, \\[10pt]
     \theta(8m_1+r_1+i,\ 8m_2+r_2+j), & \mathrm{if} & k=64+8i+j.
\end{array} \right.
\end{equation}
Thus, the renormalisation equations for the 3-point correlations can be rewritten as
\[\bs{\eta}(8m_1+r_1,\ 8m_2+r_2) = \bs{B}_{(r_1,r_2)} \ \bs{\eta}(m_1,\ m_2).  \]
We describe $\bs{B}_{(0,0)}$ and profit from the same trick as above to get $\bs{B}_{(r_1,r_2)}$. The matrix $\bs{B}_{(0,0)}\in \Mat(128,\QQ)$ can be decomposed as a sum of Kronecker products of matrices $M_{(i,j)}$ and $N_{(i,j)}$.  Direct calculation shows that the Klein four-group $K$ naturally appears in the structure of $\bs{B}_{(0,0)}$. Denote by $g(a,b)$  the standard action of $K$ on the tuple $(a,b)\in \{0,1\}^2$ and let $\tau \in K$ be as defined above.  Then, we obtain the desired decomposition

\begin{equation}
\label{eq:matrix_dim2}
    8\bs{B}_{(0,0)} = \sum_{g\in K}\begin{pmatrix}N_{g(0,0)} & N_{g(0,1)} \\[10pt] N_{g(1,0)} & N_{g(1,1)} \end{pmatrix} \otimes N_{g\tau(0,0)}+ \begin{pmatrix}M_{g(0,0)} & M_{g(0,1)} \\[10pt] M_{g(1,0)} & M_{g(1,1)} \end{pmatrix} \otimes M_{g(0,0)} .
\end{equation}

Then, the general matrix $\bs{B}_{(r_1,r_2)}$ with $r_i \in \{0,\dots,\ 7\}$ can be expressed as
\begin{multline*}
    8\bs{B}_{(r_1,r_2)} = \sum_{g\in K}\mathcal{J}^{r_1}\!\begin{pmatrix}N_{g(0,0)} & N_{g(0,1)} \\[10pt] N_{g(1,0)} & N_{g(1,1)} \end{pmatrix} \otimes \mathcal{J}^{r_2}(N_{g\tau(0,0)}) \\
    +{} \mathcal{J}^{r_1}\!\begin{pmatrix}M_{g(0,0)} & M_{g(0,1)} \\[10pt] M_{g(1,0)} & M_{g(1,1)} \end{pmatrix} \otimes \mathcal{J}^{r_2}(M_{g(0,0)}) .
\end{multline*}

The Kronecker product structure is a consequence of our choice of the vector $\bs{\eta}$ and is not surprising. The extension of the relation \eqref{eq:eta_def} for higher-order correlations is straightforward. Unfortunately, the decomposition of the matrices $\bs{B}_{(0,\dots, 0)}$, similar to \eqref{eq:matrix_dim2}, is a difficult task, mostly due to the fact that the dimension of matrix $\bs{B}_{(0,\dots, 0)}$ representing $n$-point correlations is $2\cdot8^{n-1}$. Thus in the case of 4-point correlations, we already obtain a~matrix of size $1024\times1024$. The decomposition, in this case, is still obtainable and results in 

\begin{multline*}
    16\,\mathbf{B}_{(0,0,0)} = \sum_{g\in K}\sum_{(i,j)\in \{0,1\}^2}
    \begin{pmatrix}N_{g(0,0)} & N_{g(0,1)} \\[10pt] N_{g(1,0)} & N_{g(1,1)} \end{pmatrix} \otimes N_{g(i,j)} \otimes N_{(i,j)} \\ 
    + {} \begin{pmatrix}M_{g(0,0)} & M_{g(0,1)} \\[10pt] M_{g(1,0)} & M_{g(1,1)} \end{pmatrix} \otimes M_{g(i,j)} \otimes M_{\tau(i,j)}. 
\end{multline*}
As in the 3-point correlation case, we sum over all elements of Klein four-group and all pairs $(i,j)\in \{0,1\}^2$. An increase of the dimension by 1 leads to four times more terms in the summation. On the other hand, there is still an open question about the general decomposition formula for $n$-point correlation matrix $\bs{B}_{(0, \dots, \ 0)}$, which we cannot answer at present. 

\section{4-point correlation functions}
We have already obtained an exact form of the correlation functions for 2 points and an odd number of points. The remaining cases are more complex. As discussed in the previous section, even the matrix description does not suffice to bring up a general formula for the renormalisation relations for arbitrary $n$. On the other hand, for fixed $n$, we can still study the given $n$-point correlation functions without knowing the general formulas. This section focuses on the first non-trivial higher-order correlation function, namely the 4-point one. We show the asymptotic behaviour of its sums and describe points in which the correlation function reaches the desired value. 

First, we state several facts immediately following the recursive structure of renormalisation equations and the solution of the self-consistent part.
\begin{fact}
The values of the $4$-point correlation function of the RS sequence form a~proper subset of dyadic rationals. \exend
\end{fact}

\begin{fact}
The correlation functions fulfill 
\[ |\eta^{(4)}(m_1,m_2,m_3)| \leqslant 1 \quad \mbox{and} \quad |\theta^{(4)}(m_1,m_2,m_3)| \leqslant 1  \] 
for all triples $(m_1,m_2,m_3)\in\ZZ^3$.
\exend
\end{fact}

Now, we can start our discussion of the \emph{level sets}, i.e., those subsets in $\ZZ^3$ on which the function $\eta^{(4)}$ is constant.  
\begin{prop}
\label{prop:level_set_1}
For the $4$-point correlation function of the balanced Rudin--Shapiro sequence, we have
\[ \eta^{(4)}(m_1,m_2,m_3)=1 \quad \Longleftrightarrow \quad m_1 = 0 \ \mbox{and} \ m_2=m_3, \]
up to a permutation of the indices. Moreover, $|\theta^{(4)}(m_1,m_2,m_3)| < 1$ holds for all triples $(m_1,m_2,m_3)\in\ZZ^3$.

\begin{proof}
The right to left is always true, since $w_iw_iw_{i+k}w_{i+k} =1$ holds for any $i, k\in \ZZ$. For the converse, note that the unique ergodicity of $\XX_{\mathrm{RS},2}$ implies the existence of a unique strictly positive translation invariant measure, the patch frequency measure. 
If in the sum in the definition \eqref{eq:def_corr_n} appeared any patch equals to $-1$, it would occur with a strictly positive frequency. 
Thus, the result would be strictly smaller than 1. Therefore, one needs to find a~triple of indices $(m_1,m_2,m_3)\in\ZZ^{3}$ such that $w_iw_{i+m_1}w_{i+m_2}w_{i+m_3} =1$ holds for any $i\in \ZZ$. It immediately follows that $m_j=0$  for some $j$ and therefore (since the two-point correlation does not vanish at 0 only), one has $m_{\ell}-m_k=0$ for the two remaining indices. The second claim follows immediately via the same argument.
\end{proof}
\end{prop}

We begin our further discussion on the level sets with an observation 
\[ \eta^{(4)}(1,2,3) = -\myfrac{1}{2}. \]
The renormalisation equation (with $M = m_1+m_2+m_3$)
\begin{equation}
    \label{eq:renorm_1}
    \eta^{(4)}(4m_1,4m_2,4m_3) = \myfrac{1}{2}(1+(-1)^M) \eta^{(4)}(m_1,m_2,m_3)
\end{equation}
gives for all $m\in\ZZ$ 
\[ \eta^{(4)}\bigl(4^m(1,2,3)\bigr) = -\myfrac{1}{2}. \]
Proposition \ref{prop:level_set_1} describes the 1-level set which we can use for obtaining other coordinates $(x,y,z)$ for which $\eta^{(4)}(x,y,z) = -\tfrac{1}{2}$. The renormalisation equation 
\begin{multline*}
    \eta^{(4)}(4m_1+1,4m_2+2,4m_3+3) = \myfrac{1}{4}\bigl(-(-1)^{m_2+m_3}\eta^{(4)}(m_1,m_2,m_3)\\   - (-1)^{m_1+m_2}\eta^{(4)}(m_1,m_2,m_3+1) 
      -(-1)^{m_1}\eta^{(4)}(m_1,m_2+1,m_3+1) \\ 
      + (-1)^{m_3}\eta^{(4)}(m_1+1,m_2+1,m_3+1)\bigr)
\end{multline*}
with the choice $m_1=0$, $m_2=m_3=\ell$ yields 
\begin{equation}
    \label{eq:level_set_1/2}
\eta^{(4)}  (1,4\ell+2,4\ell+3) = -\myfrac{1}{2}. 
\end{equation}
The other possible (non-trivial) choice $m_1=-1$, $m_2=m_3=\ell$ results after a translation and a permutation in 
$\eta^{(4)}  (-1,4\ell+1,4\ell+2) = \eta^{(4)}  (1,4\ell+2,4\ell+3)$.
Another renormalisation equation, namely, 
\begin{multline}
\label{eq:renorm_eq_B}
\eta^{(4)}(4m_1,4m_2+2,4m_3+2) = \myfrac{1}{2}\bigl((-1)^{m_2+m_3}
        \eta^{(4)}(m_1,m_2,m_3) \\ +(-1)^{m_1}\eta^{(4)}(m_1,m_2+1,m_3+1)\bigr) 
\end{multline}
gives with the choice $m_1=m_2=2\ell+1$, $m_3=0$ (after a permutation) 
\[\eta^{(4)}(2,8\ell+4,8\ell+6) = -\myfrac{1}{2}. \]
The other possible choices ($m_1=m_3=2\ell+1$, $m_2=0$, and $m_1=2\ell+1$, $m_2=2\ell$, $m_3=-1$) result in the same set of coordinates. 

Thus, up to now, we have found an infinite set of coordinates 
\begin{equation}
    \label{eq:level_set_1/2_2}
    \left\{ 2^m (1,4\ell+2,4\ell+3) \ : \ m\in\NN, \ \ell\in\ZZ\right\}
\end{equation}
which is a subset of the $\bigl(-\tfrac{1}{2}\bigr)$-level set. We obtained the coordinates \eqref{eq:level_set_1/2_2} as a result of studying renormalisation equations and we used those of them allowing the decomposition of the  $-\tfrac{1}{2}$ of the form $-\tfrac{1}{2} = \tfrac{1}{4}(-1-1)$ in the first case, and $-\tfrac{1}{2} = \tfrac{1}{2} (-1)$ in the second one. There are other ways to decompose $-\tfrac{1}{2}$ as a quarter of sum of four terms, for example $-\tfrac{1}{2} = \tfrac{1}{4}\bigl(-1-\tfrac{1}{2}-\tfrac{1}{2}\bigr)$. One can check that the assumptions on the $\bigl(-\tfrac{1}{2}\bigr)$-level set being equal to \eqref{eq:level_set_1/2_2} is fully consistent with the renormalisation relations in the sense, that no other coordinates can appear in the  $\bigl(-\tfrac{1}{2}\bigr)$-level set. The numerical simulations suggest that there are no other points in this level set, but a proof of this statement is still missing. 

We can proceed further with a scheme describing a way to obtain other infinite series of points belonging to the same level set. 

In order to study $\eta^{(4)}(1,2\ell_0, 2\ell_0+1)$ using the renormalisation equations, one has to treat the odd and even $\ell_0$ separately. Indeed, taking $\ell_0=2\ell_1+1$ results in \eqref{eq:level_set_1/2}. It remains to discuss the even case, i.e. $\ell_0 = 2\ell_1$. The renormalisation equation then reads
\begin{equation}
    \label{eq:renorm_eq_A}
    \eta^{(4)}(1,4\ell_1,4\ell_1+1) = \myfrac{1}{4}\bigl((2+(-1)^{\ell_1})\eta^{(4)}(0,\ell_1,\ell_1)+(-1)^{\ell_1}\eta^{(4)}(1,\ell_1,\ell_1 +1) \bigr).
\end{equation}
To solve this equation, one has to distinguish the odd and even cases again. 
The process described above can be visualised as a tree as follows.
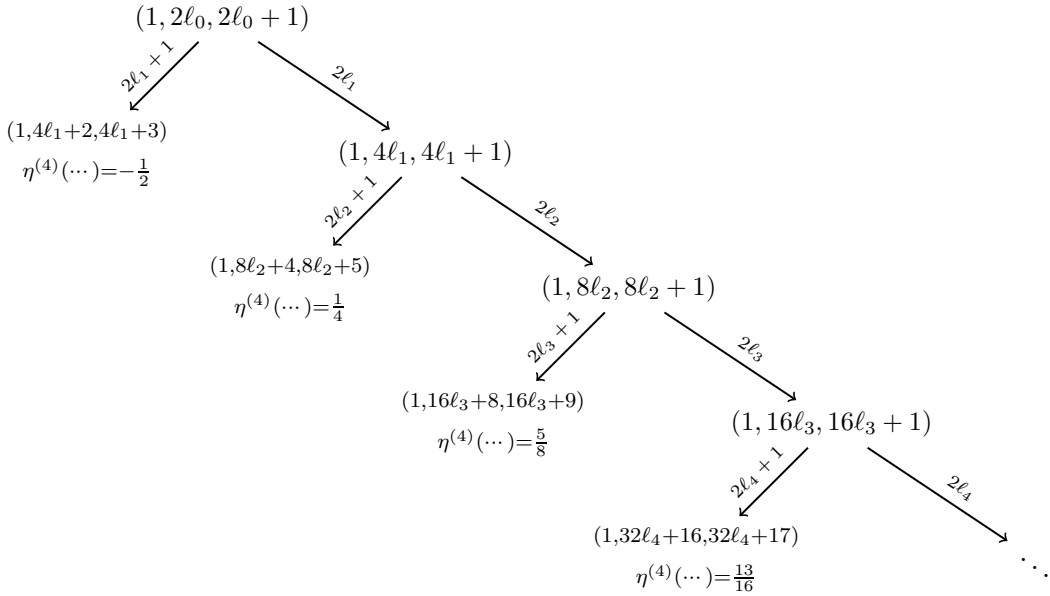
\vspace{0.5cm}
\begin{center}
\begin{tikzpicture}[scale=1.8]
    \node (0) at   (2,0) {\small$(1,2\ell_0,2\ell_0+1)$};
    \node (10) at  (1,-1) {$\substack{(1,4\ell_1+2,4\ell_1+3)\\[5pt] \eta^{(4)} (\cdots) = -\frac{1}{2}}$};
        \draw[->] (0) -- node[midway, above right, sloped, pos=1]{\tiny$2\ell_1+1$} (10);
    \node (01) at (3.5,-1) {\small$(1,4\ell_1,4\ell_1+1)$};
        \draw[->] (0) -- node[midway, above left, sloped, pos=0.75]{\tiny$2\ell_1$} (01);
    \node (20) at  (2.5,-2) {$\substack{(1,8\ell_2+4,8\ell_2+5)\\[5pt] \eta^{(4)} (\cdots) = \frac{1}{4}}$};
        \draw[->] (01) -- node[midway, above right, sloped, pos=1]{\tiny$2\ell_2+1$} (20);
    \node (02) at (5,-2) {\small$(1,8\ell_2,8\ell_2+1)$};
        \draw[->] (01) -- node[midway, above left, sloped, pos=0.75]{\tiny$2\ell_2$} (02);
    \node (30) at  (4,-3) {$\substack{(1,16\ell_3+8,16\ell_3+9)\\[5pt] \eta^{(4)} (\cdots) = \frac{5}{8}}$};
        \draw[->] (02) -- node[midway, above right, sloped, pos=1]{\tiny$2\ell_3+1$} (30);
    \node (03) at (6.5,-3) {\small$(1,16\ell_3,16\ell_3+1)$};
        \draw[->] (02) -- node[midway, above left, sloped, pos=0.75]{\tiny$2\ell_3$} (03);
    \node (40) at  (5.5,-4) {$\substack{(1,32\ell_4+16,32\ell_4+17)\\[5pt] \eta^{(4)} (\cdots) = \frac{13}{16}}$};
        \draw[->] (03) -- node[midway, above right, sloped, pos=1]{\tiny$2\ell_4+1$} (40);
    \node (04) at (8,-4) {$\ddots$};
        \draw[->] (03) -- node[midway, above left, sloped, pos=0.75]{\tiny$2\ell_4$} (04);
\end{tikzpicture}
\end{center}
The relation $\eta^{(4)}(1,2m+1,2m+2) = 0$ is often used in the calculation and covers the case $\eta^{(4)}(1,\ell,\ell+1)$ for $\ell $ odd.
Further, we can derive a formula for the values of $\eta^{(4)}$ at each level and obtain a sequence of level sets whose correlations grow towards the maximal possible value. 

\begin{prop}
\label{prop:level_set_gen1}
For the balanced Rudin--Shapiro sequence, we have 
\begin{align*}
\eta^{(4)}\bigl(1,\, 2^{n}(2\ell+1),\,2^{n}(2\ell+1)+1 \bigr) \, &{}= \, 1-\myfrac{3}{2^n}  \\
\eta^{(4)}(1,\,2\ell+1,\,2\ell+2) &{}= 0
\end{align*}
for each $n\in\NN$ and for $\ell \in \ZZ$. 
\begin{proof}
The second identity can be obtained by separately treating odd and even $\ell$. We show the first claim by induction on $n$. For $n=1$ we have already proved that $\eta^{(4)}(1,4\ell+2,4\ell+3) \ =\ -\tfrac{1}{2} \ =\ 1 - \tfrac{3}{2}$. 
For $n=2$ we get $\eta^{(4)}(1,8\ell+4,8\ell+5)$ and using  \eqref{eq:renorm_eq_A} one has
$\eta^{(4)}\bigl(1,4(2\ell+1),4(2\ell+1)+1 \bigr) = \tfrac{1}{4} = 1-\tfrac{3}{4}. $
Now, suppose that the claim holds for any $n'\leqslant n-1$ and recall the equation \eqref{eq:renorm_eq_A}. Then, for any $n\geqslant 3$, one obtains 
\begin{align*}
    \eta^{(4)}&\bigl(1,\, 2^{n}(2\ell+1),\,2^{n}(2\ell+1)+1 \bigr) = \eta^{(4)}\bigl(1,\, 4(2^{n-2}(2\ell+1)),\, 4(2^{n-2}(2\ell+1))+1 \bigr) \\[5pt]
    &= \myfrac{1}{4}\bigr(3\,\eta^{(4)}(0,\, 2^{n-2}(2\ell+1),\, 2^{n-2}(2\ell+1)) + \eta^{(4)} (1,\, 2^{n-2}(2\ell+1),\, 2^{n-2}(2\ell+1)+1 ) \bigl) \\[5pt]
    &=\myfrac{1}{4}\bigl(3+1-\myfrac{3}{2^{n-2}} \bigr) = 1-\myfrac{3}{2^n},
\end{align*}
using $\eta^{(4)} (1,\, 2^{n-2}(2\ell+1),\, 2^{n-2}(2\ell+1)+1 ) = 1-\frac{3}{2^{n-2}}$. 
\end{proof}
\end{prop}

We can profit from this result and extend the current level sets and ``double" their cardinality.  
\begin{prop}
\label{prop:level_set_gen2}
For the balanced Rudin--Shapiro sequence, we have 
\[\eta^{(4)}\bigl(2,\, 2^{n+1}(2\ell+1),\,2^{n+1}(2\ell+1)+2 \bigr) \, = \, 1-\myfrac{3}{2^n}  \]
for each $n\in\NN$ and $\ell \in \ZZ$. 
\end{prop}
\pushQED{\qed}
\noindent \textit{Proof.}
We already know that $\eta^{(4)}(2,8\ell+4,8\ell+6) = -\myfrac{1}{2}$. For $n\geqslant 2$ one gets with help of \eqref{eq:renorm_eq_B} and Proposition \ref{prop:level_set_gen1} the result. \begin{align*}
    \eta^{(4)}&\bigl(2, 2^{n+1}(2\ell+1),2^{n+1}(2\ell+1)+2 \bigr) \\
    &=\myfrac{1}{2}\bigl( \eta^{(4)}\bigl(0, 2^{n-1}(2\ell+1),2^{n-1}(2\ell+1) \bigr) + \eta^{(4)}\bigl(1, 2^{n-1}(2\ell+1),2^{n-1}(2\ell+1)+1 \bigr) \bigr) \\
    &= \myfrac{1}{2} \bigl(1+1- \myfrac{3}{2^{n-1}}\bigr) = 1-\myfrac{3}{2^{n}}. \qedhere 
\end{align*} 
\popQED

\noindent Observe that the sum of coordinates  $1+4(2\ell+1)+4(2\ell+1)+1$ is always even and so does $2+2^{n+1}(2\ell+1)+2^{n+1}(2\ell+1)+2$. Thus, we can recall the renormalisation equation \eqref{eq:renorm_1} and get the final description of certain sets in the positive octant where the function $\eta^{(4)}$ is constant. We can further extend these sets to all of $\ZZ^3$ via Proposition \ref{prop:symmetries}.

\begin{theorem}
\label{thm:level_set}
The $4$-point correlation function $\eta^{(4)}$ of the Rudin--Shapiro sequence with balanced weights, for every $n\geqslant 1$, is constant on the set
\[ \mathcal{C}_n \eqdef \left\{ 2^{m} \bigl(1,\, 2^{n}(2\ell+1),\,2^{n}(2\ell+1)+1 \bigr) \ \mbox{and all permutations}  \ : \ m\geqslant 0, \ \ell\in\ZZ \right\}\]
and reaches the value 
\[\pushQED{\qed}\eta^{(4)}{\Big|}_{\mathcal{C}_n} \ = \ 1- \myfrac{3}{2^{n}}. \qedhere \popQED \]
\end{theorem}
Similar considerations as above lead to a description of the function $\theta^{(4)}$ evaluated at points from the set $\mathcal{C}_n$. Note that if we want to extend the results to $\ZZ^3$ in this case, the translation symmetry may add an additional minus factor (as stated in Proposition \ref{prop:symmetries}). 
\begin{prop}
For the function $\theta^{(4)}$ of the Rudin--Shapiro sequence with balanced weights, one has, for all $\ell \in \ZZ$, $n\in\NN$ and $m\in \NN_0$,
\begin{align*}
\pushQED{\qed}
\theta^{(4)}\bigl(2^{m}(1,\, 2^{n}(2\ell+1),\,2^{n}(2\ell+1)+1 ) \bigr) &{}=  
     \left( \frac{3}{2^{n}}-\delta_{_{n,1}}\right)\delta_{_{0,m}}, \\[5pt]
\theta^{(4)}\bigl(2^{m}(1,\, 2\ell+1,\,2\ell+2 ) \bigr)&{} =  0.   \qedhere
\end{align*}
\popQED
\end{prop}

Of course, the result also holds for all permutations of the coordinates, but we do not repeat this in the upcoming propositions. 

It turns out that the strategy described above can be applied to various ``starting" vectors (i.e., different to $(1,\ell,\ell+1)$), and one gets the description of the correlation functions at different infinite subsets of $\ZZ^3$. 
The proofs of the following propositions are technical and follow the above scheme  (using suitable renormalisation equations), and profit from the results of Theorem \ref{thm:level_set}. Therefore, we decided to omit them. 
\begin{prop}[Starting vector $(1,\ell,\ell+2)$]
\label{prop:vectors1}
For the functions $\eta^{(4)}$, $\theta^{(4)}$ of the Rudin--Shapiro sequence with balanced weights one has for all $\ell \in \ZZ$
\begin{align*}
\pushQED{\qed}
\eta^{(4)}\bigl(2^{m}(1,\, 2^{n}(2\ell+1),\,2^{n}(2\ell+1)+2 ) \bigr) {}&=  \left\{\begin{array}{rl}
    0, & \mbox{ if } m=0 \mbox{ and } n=1, \\[5pt]
     \frac{1}{4}, & \mbox{ if } m=0 \mbox{ and } n=2, \\[5pt]
     -\frac{3}{2^{n}}, & \mbox{ if } m=0 \mbox{ and } n\geqslant 3,\\[5pt]
     0, & \mbox{ if } m\geqslant 1.
\end{array} \right. \\[10pt]
\theta^{(4)}\bigl(2^{m}(1,\, 2^{n}(2\ell+1),\,2^{n}(2\ell+1)+2 ) \bigr) {}&=  \left\{\begin{array}{rl}
    0, & \mbox{ if } m\in\{0,1\} \mbox{ and } n=1, \\[5pt]
     -(-1)^m \frac{1}{4}, & \mbox{ if } m\in\{0,1\} \mbox{ and } n=2, \\[5pt]
     (-1)^m\frac{3}{2^{n}}, & \mbox{ if } m\in\{0,1\} \mbox{ and } n\geqslant 3,\\[5pt]
     0, & \mbox{ if } m\geqslant 2.
\end{array} \right. \\[10pt]
    \eta^{(4)}\bigl(2^{m}(1,\, 2^{n}(2\ell+1)-1,\,2^{n}(2\ell+1)+1 ) \bigr) {}&=  \left\{\begin{array}{rl}
    0, & \mbox{ if } m=0 \mbox{ and } n=1, \\[5pt]
     -\frac{1}{4}, & \mbox{ if } m=0 \mbox{ and } n=2, \\[5pt]
     \frac{3}{2^{n}}, & \mbox{ if } m=0 \mbox{ and } n\geqslant 3,\\[5pt]
     0, & \mbox{ if } m\geqslant 1.
\end{array} \right. \\[10pt]
\theta^{(4)}\bigl(2^{m}(1,\, 2^{n}(2\ell+1)-1,\,2^{n}(2\ell+1)+1 ) \bigr) {}&=  \left\{\begin{array}{rl}
    0, & \mbox{ if } m\in\{0,1\} \mbox{ and } n=1, \\[5pt]
     (-1)^m \frac{1}{4}, & \mbox{ if } m\in\{0,1\} \mbox{ and } n=2, \\[5pt]
     -(-1)^m\frac{3}{2^{n}}, & \mbox{ if } m\in\{0,1\} \mbox{ and } n\geqslant 3,\\[5pt]
     0, & \mbox{ if } m\geqslant 2. \qquad \qquad \qquad \quad \qedhere
\end{array} \right. 
\end{align*}
\popQED
\end{prop}

\begin{prop}[Starting vector $(2,\ell,\ell+1)$]
\label{prop:vectors2}
For the functions $\eta^{(4)}$, $\theta^{(4)}$ of the Rudin--Shapiro sequence with balanced weights one has for all $\ell \in \ZZ$
\begin{align*}
\eta^{(4)}\bigl(2^{m}(2,\, 2^{n}(2\ell+1),\,2^{n}(2\ell+1)+1 ) \bigr) {}&=  \left\{\begin{array}{rl}
     0, & \mbox{ if } m=0 \mbox{ and } n=1, \\[5pt]
     \frac{1}{4}, & \mbox{ if } m=0 \mbox{ and } n= 2,\\[5pt]
     -\frac{3}{2^{n}}, & \mbox{ if } m=0 \mbox{ and } n\geqslant 3,\\[5pt]
     0, & \mbox{ if } m\geqslant 1,
\end{array} \right. \\[10pt]
\theta^{(4)}\bigl(2^{m}(2,\, 2^{n}(2\ell+1),\,2^{n}(2\ell+1)+1 ) \bigr) {}&=  \left\{\begin{array}{rl}
     0, & \mbox{ if } m\in\{0,1\} \mbox{ and } n=1, \\[5pt]
     -\frac{1}{4}, & \mbox{ if } m\in\{0,1\}  \mbox{ and } n= 2,\\[5pt]
     \frac{3}{2^{n}}, & \mbox{ if } m\in\{0,1\}  \mbox{ and } n\geqslant 3,\\[5pt]
     0, & \mbox{ if } m\geqslant 2,
\end{array} \right. \end{align*}
\begin{align*} \pushQED{\qed}
\eta^{(4)}\bigl(2^{m}(2,\, 2^{n}(2\ell+1)+1,\,2^{n}(2\ell+1)+2 ) \bigr) {}&=  \left\{\begin{array}{rl}
     0, & \mbox{ if } m=0 \mbox{ and } n=1, \\[5pt]
     -\frac{1}{4}, & \mbox{ if } m=0 \mbox{ and } n= 2,\\[5pt]
     \frac{3}{2^{n}}, & \mbox{ if } m=0 \mbox{ and } n\geqslant 3,\\[5pt]
     0, & \mbox{ if } m\geqslant 1,
\end{array} \right. \\[10pt]
\theta^{(4)}\bigl(2^{m}(2,\, 2^{n}(2\ell+1)+1,\,2^{n}(2\ell+1)+2 ) \bigr) {}&=  \left\{\begin{array}{rl}
     0, & \mbox{ if } m\in\{0,1\} \mbox{ and } n=1, \\[5pt]
     -(-1)^{m}\frac{1}{4}, & \mbox{ if } m\in\{0,1\}  \mbox{ and } n= 2,\\[5pt]
     (-1)^{m}\frac{3}{2^{n}}, & \mbox{ if } m\in\{0,1\}  \mbox{ and } n\geqslant 3,\\[5pt]
     0, & \mbox{ if } m\geqslant 2. \qquad \qquad \qquad \quad \qedhere
\end{array} \right. 
\end{align*}
\popQED
\end{prop}

\begin{prop}[Starting vector $(2,\ell,\ell+2)$]
\label{prop:vectors3}
For the functions $\eta^{(4)}$, $\theta^{(4)}$ of the Rudin--Shapiro sequence with balanced weights one has for all $\ell \in \ZZ$ and for all $m,n\in\NN_{0}$

\begin{align*}\pushQED{\qed}
    \eta^{(4)}\bigl(2^{m}(2,\, 2^{n+2}(2\ell+1)+1,\,2^{n+2}(2\ell+1)+3 ) \bigr) {}&= -\myfrac{1}{2} + \myfrac{3}{2^{n+2}}, \\[5pt]
    \theta^{(4)}\bigl(2^{m}(2,\, 2^{n+2}(2\ell+1)+1,\,2^{n+2}(2\ell+1)+3 ) \bigr) {}&= \delta_{_{m,0}}\bigl(-\myfrac{1}{2} + \myfrac{3}{2^{n+2}}  \bigr), \\[5pt]
    \eta^{(4)}\bigl(2^{m}(2,\, 2^{n+2}(2\ell+1)-1,\,2^{n+2}(2\ell+1)+1 ) \bigr) {}&= -\myfrac{1}{2} + \myfrac{3}{2^{n+2}}, \\[5pt]
    \theta^{(4)}\bigl(2^{m}(2,\, 2^{n+2}(2\ell+1)-1,\,2^{n+2}(2\ell+1)+1 ) \bigr) {}&= \delta_{_{m,0}}\bigl(-\myfrac{1}{2} + \myfrac{3}{2^{n+2}}  \bigr). \qedhere 
\end{align*}
\popQED
\end{prop}

Similarly, one can continue this procedure and generate infinite series where the correlations function remains constant. Theorem \ref{thm:level_set} and Propositions \ref{prop:vectors1}, \ref{prop:vectors2}, \ref{prop:vectors3} provide infinitely many infinite arithmetic progressions in $\ZZ^3$ along which is $\eta^{(4)}$ constant and non-zero.  Therefore, if we move to the averages of the correlation functions, one might expect that the averages over the cube $\{0,1,\dots N-1\}^3$ for any $N$ cannot be arbitrarily small. In what follows, we prove the opposite, namely, that the mean of the distances $|\eta^{(4)}|$ vanishes asymptotically, and we conclude that the average vanishes asymptotically as well.

\begin{prop}
For the $4$-point Rudin--Shapiro correlation functions, one has
\begin{align*}
\lim_{N\to \infty} \myfrac{1}{N^3} \sum_{0\leqslant m_i\leqslant N-1} \bigl|\eta^{(4)}(m_1,m_2,m_3) \bigr| &{}= 0, \\
\lim_{N\to \infty} \myfrac{1}{N^3} \sum_{0\leqslant m_i\leqslant N-1} \bigl|\theta^{(4)}(m_1,m_2,m_3) \bigr| &{}= 0. 
\end{align*}

\begin{proof}
Denote by 
\[ \mathbf{\Sigma}(N) = \myfrac{1}{N^3} \sum_{0\leqslant m_i\leqslant N-1} \bigl|\eta^{(4)}(m_1,m_2,m_3) \bigr|, \quad \mbox{and} \quad
    \mathbf{\Theta}(N) = \myfrac{1}{N^3} \sum_{0\leqslant m_i\leqslant N-1} \bigl|\theta^{(4)}(m_1,m_2,m_3) \bigr|. \]
We provide the following calculations for $\mathbf{\Sigma}$. The estimates for $\mathbf{\Theta}$ are analogous.

First, observe that since all correlations are smaller than or equal to 1, one has
\begin{align*}
    \mathbf{\Sigma}(4N+1) {}&= \myfrac{1}{(4N+1)^3}\sum_{0\leqslant m_i\leqslant 4N} \bigl|\eta^{(4)}(m_1,m_2,m_3) \bigr|\\
    {}&= \myfrac{1}{(4N+1)^3}\sum_{0\leqslant m_i\leqslant 4N-1} \bigl|\eta^{(4)}(m_1,m_2,m_3) \bigr| + \myfrac{1}{(4N+1)^3}\underbrace{ \sum_{\substack{m_i =4N \\[3pt] i \in\{1,2,3\}  }} \bigl|\eta^{(4)}(m_1,m_2,m_3) \bigr|}_{\leqslant 3(4N)^2 + 1} \\
    {}&= \myfrac{(4N)^3}{(4N+1)^3}\ \mathbf{\Sigma}(4N) + O(N^{-1})
\end{align*}
as $N \to \infty$. Analogously, similar relations hold for $ \mathbf{\Sigma}(4N+2)$ and $\mathbf{\Sigma}(4N+3)$. Hence, if we show that $\mathbf{\Sigma}(4N) \rightarrow 0$ as $N \to \infty$, our claim follows. 
Thus, consider
\begin{align*}
    \mathbf{\Sigma}(4N) {}& = \myfrac{1}{(4N)^3}\sum_{0\leqslant m_i\leqslant 4N-1} \bigl|\eta^{(4)}(m_1,m_2,m_3) \bigr|\\[2mm]
    {}& = \myfrac{1}{(4N)^3}\sum_{0\leqslant m_i\leqslant N-1} \sum_{r_i \in \{ 0,1,2,3\}} \bigl|\eta^{(4)}(4m_1+r_1,4m_2+r_2,4m_3+r_3) \bigr|\\
    {}& \leqslant \myfrac{1}{4(4N)^3}\left(\sum_{0\leqslant m_i\leqslant N-1} \left( 128\,\bigl|\eta^{(4)}(m_1,m_2,m_3) \bigr|+112\,\bigl|\theta^{(4)}(m_1,m_2,m_3) \bigr|\right) + O(N^2)\right) \\[10pt]
    {}& =\myfrac{1}{4^4}\bigl( 128\,\mathbf{\Sigma}(N) + 112\,\mathbf{\Theta}(N)\bigr) + O(N^{-1}).
\end{align*}
In the third row, we inserted the renormalisation equations and  used the triangle inequality for each of them together with
\[ \sum_{0\leqslant m_i\leqslant N-1} \eta^{(n)}(m_1+r_1,\dots,m_{n-1}+r_{n-1}) = \sum_{0\leqslant m_i\leqslant N-1} \eta^{(n)}(m_1,\dots,m_{n-1}) + O(N^{n-2}), \]
which holds \footnote{It suffices because the only term of this form appears on the RHS of the renormalisation equations.} for any $(r_1,\dots,r_{n-1}) \in \{0,1\}^{n-1}$. Analogous estimates hold for $\theta^{(n)}$ as well. 
For the summatory function $\mathbf{\Theta}$, we get the following relation
\[\mathbf{\Theta}(4N) \leqslant \myfrac{1}{4^4}\bigl( 112\,\mathbf{\Sigma}(N) + 128\,\mathbf{\Theta}(N)\bigr) + O(N^{-1}). \]
Combining both equations, one gets
\[\mathbf{\Sigma}(4N) + \mathbf{\Theta}(4N) \leqslant \myfrac{15}{16}\bigl(\mathbf{\Sigma}(N) + \mathbf{\Theta}(N) \bigr) + O(N^{-1}). \]
This inequality implies $\lim_{N\to+\infty} \mathbf{\Sigma}(4N) + \mathbf{\Theta}(4N) =0 $. Since $\mathbf{\Sigma}(4N)$ and $\mathbf{\Theta}(4N)$ are positive, we get  the desired, namely,
$\lim_{N\to+\infty} \mathbf{\Sigma}(4N) =0 $ and $\lim_{N\to+\infty} \mathbf{\Theta}(4N) =0 $.
\end{proof}
\end{prop}

This convergence has an immediate consequence: the triangle inequality gives the desired result on the asymptotic behaviour of the mean of the coefficients. Moreover, since all coefficients are in modulus smaller than or equal to one, we obtain the asymptotically vanishing means for arbitrary powers $\alpha \geqslant 1$ of the correlation functions. 

\begin{coro}
For the $4$-point Rudin--Shapiro correlation functions and any $\alpha \geqslant 1$, one has
\begin{align*}\pushQED{\qed}
\lim_{N\to \infty} \myfrac{1}{N^3} \sum_{0\leqslant m_i\leqslant N-1} \eta^{(4)}(m_1,m_2,m_3)^{\alpha}  &{}= 0, \\
\lim_{N\to \infty} \myfrac{1}{N^3} \sum_{0\leqslant m_i\leqslant N-1} \theta^{(4)}(m_1,m_2,m_3)^{\alpha} &{}= 0. \qedhere
\end{align*}
\popQED
\end{coro}

\begin{remark}
Using the renormalisation equations without the absolute value and triangle inequality, one could improve the estimates for the sums of correlation functions (and not their absolute values). We leave this part to interested readers. 
\exend
\end{remark}

Even though we showed that the average of 4-point correlations is zero, we can still recognise a difference from a random structure. The presence of infinitely many infinitely long arithmetic progressions with a constant non-zero value of $\eta^{(4)}$ suggests a presence of a certain long-range ordering. We show an explicit example of a structure in the original RS sequence detected by the 4-point correlation function. 

The $\bigl(-\tfrac{1}{2}\bigr)$-level set contains as its subset vectors of the form $(1,4k+2,4k+3)$. It is worth evaluating the 4-point correlation function at the point $(1,k,k+1)$. Proposition \ref{prop:level_set_gen1} provides its description, but we wish to have the renormalisation equations for these particular vectors at hand. They form a~closed set of equations, namely
\begin{align*}
    \eta^{(4)}(1,4m,4m+1) & = \myfrac{1}{4}\bigl(2+(-1)^m(1+\eta^{(4)}(1,m,m+1))\bigr),\\
    \eta^{(4)}(1,4m+1,4m+2) & =0, \\
    \eta^{(4)}(1,4m+2,4m+3) & = - \myfrac{1}{2}, \\
    \eta^{(4)}(1,4m+3,4m+4) & =0. 
\end{align*}

The 4-point correlation function evaluated at the point $(1,k,k+1)$ effectively measures the correlation between two doubles of two consecutive points at a distance of $k$. Nevertheless, these doubles are nothing but elements of a fixed point of the induced two-letter substitution $\varrho^{}_{2}$ defined in \eqref{eq:subst_induced}, which can also be studied in terms of its correlation functions. The corresponding set of renormalisation equations reads 
\begin{align*}
    \eta^{}_2(4k) & = \myfrac{1}{4}\bigl(2+(-1)^k(1+\eta^{}_2(k))\bigr),\\
    \eta^{}_2(4k+1) & =0, \\
    \eta^{}_2(4k+2) & = - \myfrac{1}{2}, \\
    \eta^{}_2(4k+3) & =0. 
\end{align*}
To derive them, we used the relations from Remark \eqref{eq:derived-description} and the fact
that one has
\[\lim_{N\to \infty} \myfrac{1}{N} \sum_{i=k}^{k+N-1}w_i = \lim_{N\to \infty} \myfrac{1}{N} \sum_{i=k}^{k+N-1}(-1)^i w_i = 0 \qquad \mbox{for all $k\in\ZZ$}. \]
Comparing this set of equations for $\eta^{(4)}(1,m,m+1)$ and $\eta^{}_2 (m)$ together with the corresponding initial conditions, one can see that the functions $\eta^{(4)}(1,m,m+1)$ and $\eta^{}_2(m)$ coincide. This observation illustrates that the high-order correlation functions at certain points can be understood as ordinary correlation functions for patches in the original sequence. 

In summary, we have exploited the structure of higher-order correlation functions for the binary RS sequence. Using the renormalisation approach, we have shown that all odd-point correlations vanish, and for arbitrary even $n>2$, we have found a non-zero point where the correlation differs from zero. On the other hand, we have proved that the average of these coefficients as well as the average of their distances equal zero. These results provide a better understanding of the statistical differences between the RS sequence and a random binary one. Further, we have given a detailed description of 4-point correlations and shown that they contain many arithmetic structures that can (and should) be studied further. One can ask what further symmetries can be found within them and how they can help us complete the description of 4-point correlations.


\section*{Acknowledgements}
I would like to thank Michael Baake for encouraging me to start studying the correlation functions of RS sequences and for several discussions. I also express my gratitude to the two anonymous referees for valuable suggestions that helped to improve the manuscript. This work was supported by the German Research Foundation (DFG) within the CRC 1283/2 (2021 - 317210226) at Bielefeld University.

\newpage
\section*{Appendix - Renormalisation equations for 3-point correlation functions}
In this section, we omit the upper index and write $\eta, \theta$ instead of $\eta^{(n)},\theta^{(n)}$. The number of points should be clear from the context. We include the only necessary equations. 
{\footnotesize
\begin{align*}
    \eta(4m_1,4m_2) ={}& \myfrac{1}{2}\,\eta(m_1,m_2),\\
    \eta(4m_1,4m_2\!+\!1) ={}& \myfrac{1}{4} \left(\eta(m_1,m_2) +(-1)^{m_2}(1-(-1)^{m_1})\theta(m_1,m_2) + (-1)^{m_1}\eta(m_1,m_2\!+\!1) \right),\\
    \eta(4m_1,4m_2\!+\!2) ={}& \myfrac{1}{2}(-1)^{m_1}\eta(m_1,m_2+1),\\
    \eta(4m_1,4m_2\!+\!3) ={}& \myfrac{1}{4} \left(-(-1)^{m_2}\theta(m_1,m_2) + (1+(-1)^{m_1})\eta(m_1,m_2\!+\!1)-(-1)^{m_1+m_2}\theta(m_1,m_2\!+\!1) \right),\\
    \eta(4m_1\!+\!1,4m_2\!+\!1) ={}& \myfrac{1}{4} \left((1+(-1)^{m_1+m_2})\eta(m_1,m_2) + (-1)^{m_1+m_2} \theta(m_1,m_2) - \theta(m_1\!+\!1,m_2\!+\!1) \right),\\
    \begin{split}
        \eta(4m_1\!+\!1,4m_2\!+\!2) ={}& \myfrac{1}{4} \left(-(-1)^{m_1+m_2}\eta(m_1,m_2) + (-1)^{m_2} \theta(m_1,m_2) \right. \\ & \hspace{5.2cm} \left. -(-1)^{m_1} \eta(m_1,m_2\!+\!1) -\theta(m_1\!+\!1,m_2\!+\!1) \right), 
    \end{split}\\
    \begin{split}
        \eta(4m_1\!+\!1,4m_2\!+\!3) ={}& \myfrac{1}{4} \left(-(-1)^{m_2} \theta(m_1,m_2) +(-1)^{m_1} \theta(m_1,m_2\!+\!1) \right. \\ & \hspace{4cm} \left. -(-1)^{m_1} \eta(m_1,m_2\!+\!1) +(-1)^{m_2}\eta(m_1\!+\!1,m_2\!+\!1) \right), 
    \end{split}\\
    \eta(4m_1\!+\!2,4m_2\!+\!2) ={}& \myfrac{1}{2}(-1)^{m_1+m_2}\eta(m_1,m_2),\\
    \begin{split}
        \eta(4m_1\!+\!2,4m_2\!+\!3) ={}& \myfrac{1}{4} \left(-(-1)^{m_1+m_2} \eta(m_1,m_2) -(-1)^{m_1} \theta(m_1,m_2\!+\!1) \right. \\ & \hspace{4cm} \left. + (-1)^{m_2}\eta(m_1\!+\!1,m_2\!+\!1) +\theta(m_1\!+\!1,m_2\!+\!1) \right), 
    \end{split}\\
    \eta(4m_1\!+\!3,4m_2\!+\!3) ={}& \myfrac{1}{4} \left((-1)^{m_1+m_2}\eta(m_1,m_2) + \eta(m_1\!+\!1,m_2\!+\!1)  +(1-(-1)^{m_1+m_2}) \theta(m_1\!+\!1,m_2\!+\!1) \right),\\
    \theta(4m_1,4m_2) ={}& \myfrac{1}{2}(-1)^{m_1+m_2}\theta(m_1,m_2),\\
    \theta(4m_1,4m_2\!+\!1) ={}& \myfrac{1}{4} \left(\eta(m_1,m_2) -(-1)^{m_2}(1+(-1)^{m_1})\theta(m_1,m_2) - (-1)^{m_1}\eta(m_1,m_2\!+\!1) \right),\\
    \theta(4m_1,4m_2\!+\!2) ={}& \myfrac{1}{2}(-1)^{m_2}\theta(m_1,m_2),\\
    \theta(4m_1,4m_2\!+\!3) ={}& \myfrac{1}{4} \left(-(-1)^{m_2}\theta(m_1,m_2) + ((-1)^{m_1}-1)\eta(m_1,m_2\!+\!1)+(-1)^{m_1+m_2}\theta(m_1,m_2\!+\!1) \right),\\
    \theta(4m_1\!+\!1,4m_2\!+\!1) ={}& \myfrac{1}{4} \left((1-(-1)^{m_1+m_2})\eta(m_1,m_2) + (-1)^{m_1+m_2} \theta(m_1,m_2) + \theta(m_1\!+\!1,m_2\!+\!1) \right),\\
    \begin{split}
        \theta(4m_1\!+\!1,4m_2\!+\!2) ={}& \myfrac{1}{4} \left((-1)^{m_1+m_2}\eta(m_1,m_2) + (-1)^{m_2} \theta(m_1,m_2) \right. \\ & \hspace{5.2cm} \left. -(-1)^{m_1} \eta(m_1,m_2\!+\!1) +\theta(m_1\!+\!1,m_2\!+\!1) \right), 
    \end{split}\\
    \begin{split}
        \theta(4m_1\!+\!1,4m_2\!+\!3) ={}& \myfrac{1}{4} \left(-(-1)^{m_2} \theta(m_1,m_2) -(-1)^{m_1} \theta(m_1,m_2\!+\!1) \right. \\ & \hspace{4cm} \left. -(-1)^{m_1} \eta(m_1,m_2\!+\!1) -(-1)^{m_2}\eta(m_1\!+\!1,m_2\!+\!1) \right), 
    \end{split}\\
    \theta(4m_1\!+\!2,4m_2\!+\!2) ={}& \myfrac{1}{2}\theta(m_1\!+\!1,m_2\!+\!1),\\
    \begin{split}
        \theta(4m_1\!+\!2,4m_2\!+\!3) ={}& \myfrac{1}{4} \left(-(-1)^{m_1+m_2} \eta(m_1,m_2) +(-1)^{m_1} \theta(m_1,m_2\!+\!1) \right. \\ & \hspace{4cm} \left. - (-1)^{m_2}\eta(m_1\!+\!1,m_2\!+\!1) +\theta(m_1\!+\!1,m_2\!+\!1) \right), 
    \end{split}\\
    \theta(4m_1\!+\!3,4m_2\!+\!3) ={}& \myfrac{1}{4} \left((-1)^{m_1+m_2}\eta(m_1,m_2) - \eta(m_1\!+\!1,m_2\!+\!1)  +(1+(-1)^{m_1+m_2}) \theta(m_1\!+\!1,m_2\!+\!1) \right).\\
\end{align*}
}

\section*{Appendix - Renormalisation equations for 4-point correlation functions}
We denote by $M$ the sum of all indices, i.e. $M=m_1+m_2+m_3$
{\footnotesize
\begin{align*}
    \eta(4m_1,4m_2,4m_3) ={}& \myfrac{1}{2}(1\!+\!(-1)^M) \eta(m_1,m_2,m_3),\\
    \eta(4m_1,4m_2,4m_3\!+\!1) ={}& \myfrac{1}{4}\bigl((1\!-\!(-1)^M) \eta(m_1,m_2,m_3)  
        \!+\!(-1)^{m_3}\theta(m_1,m_2,m_3)\!-\!(-1)^{m_1+m_2}\theta(m_1,m_2,m_3\!+\!1)\bigr),\\
    \eta(4m_1,4m_2,4m_3\!+\!2) ={}& 0,\\
    \eta(4m_1,4m_2,4m_3\!+\!3) ={}& \myfrac{1}{4}\bigl(-\!(-1)^{m_3}\theta(m_1,m_2,m_3)\!+\!(1\!+\!(-1)^M)
        \eta(m_1,m_2,m_3\!+\!1) \!+\!(-1)^{m_1+m_2}\theta(m_1,m_2,m_3\!+\!1)\bigr),\\
    \eta(4m_1,4m_2\!+\!1,4m_3\!+\!1) ={}& \myfrac{1}{4}\bigl((1\!+\!(-1)^{m_2+m_3}\!+\!(-1)^M)
        \eta(m_1,m_2,m_3) \!+\!(-1)^{m_1}\eta(m_1,m_2\!+\!1,m_3\!+\!1)\bigr),\\
    \begin{split}
    \eta(4m_1,4m_2\!+\!1,4m_3\!+\!2) ={}& \myfrac{1}{4}\left(-\!(-1)^{m_2+m_3}\eta(m_1,m_2,m_3)\!+\! (-1)^{m_3}\theta(m_1,m_2,m_3) \right. \\ 
    & \hspace{3cm} \left. -\!(-1)^{m_1+m_2}\theta(m_1,m_2,m_3\!+\!1) \!+\!(-1)^{m_1}\eta(m_1,m_2\!+\!1,m_3\!+\!1)\right),
    \end{split}\\
    \begin{split}
    \eta(4m_1,4m_2\!+\!1,4m_3\!+\!3) ={}& \myfrac{1}{4}\left( -\!(-1)^{m_3}\theta(m_1,m_2,m_3)\!+\!(-1)^{m_2}(1\!-\!(-1)^{m_1})\theta(m_1,m_2,m_3\!+\!1) \right. \\
    & \hspace{3cm} \left. -\!(-1)^{m_1+m_3}\theta(m_1,m_2\!+\!1,m_3\!+\!1)\right),
    \end{split}\\
    \eta(4m_1,4m_2\!+\!2,4m_3\!+\!2) ={}& \myfrac{1}{2}\bigl((-1)^{m_2+m_3}
        \eta(m_1,m_2,m_3) \!+\!(-1)^{m_1}\eta(m_1,m_2\!+\!1,m_3\!+\!1)\bigr),\\
    \begin{split}
    \eta(4m_1,4m_2\!+\!2,4m_3\!+\!3) ={}& \myfrac{1}{4}\left(-\!(-1)^{m_2+m_3}\eta(m_1,m_2,m_3)\!-\! (-1)^{m_2}\theta(m_1,m_2,m_3\!+\!1) \right. \\ 
    & \hspace{3cm} \left. +\!(-1)^{m_1}\eta(m_1,m_2\!+\!1,m_3\!+\!1) \!-\!(-1)^{m_1+m_3}\theta(m_1,m_2\!+\!1,m_3\!+\!1)\right),
    \end{split}\\
    \eta(4m_1,4m_2\!+\!3,4m_3\!+\!3) ={}& \myfrac{1}{4}\bigl((-1)^{m_2+m_3}
        \eta(m_1,m_2,m_3) \!+\!(1\!+\!(-1)^{m_1}\!+\!(-1)^M)\eta(m_1,m_2\!+\!1,m_3\!+\!1)\bigr),\\
    \eta(4m_1\!+\!1,4m_2\!+\!1,4m_3\!+\!1) ={}& \myfrac{1}{4}\bigl((1\!-\!(-1)^M) \eta(m_1,m_2,m_3)  
        \!+\!(-1)^{M}\theta(m_1,m_2,m_3)\!-\theta(m_1\!+\!1,m_2\!+\!1,m_3\!+\!1)\bigr),\\
    \begin{split}
    \eta(4m_1\!+\!1,4m_2\!+\!1,4m_3\!+\!2) ={}& \myfrac{1}{4}\big( ((-1)^{m_3}\!-\!(-1)^M)\theta(m_1,m_2,m_3)\!+\!(-1)^{m_1+m_2}\theta(m_1,m_2,m_3\!+\!1)  \\
    & \hspace{3cm} \left. - \theta(m_1\!+\!1,m_2\!+\!1,m_3\!+\!1)\right),
    \end{split}\\
    \begin{split}
    \eta(4m_1\!+\!1,4m_2\!+\!1,4m_3\!+\!3) ={}& \myfrac{1}{4}\left(-\!(-1)^{m_3}\theta(m_1,m_2,m_3)\!+\! (-1)^{m_1+m_2}\eta(m_1,m_2,m_3\!+\!1) \right. \\ 
    & \hspace{3cm} \left. +\!(-1)^{m_1+m_2}\theta(m_1,m_2,m_3\!+\!1) \!+\!(-1)^{m_3}\eta(m_1\!+\!1,m_2\!+\!1,m_3\!+\!1)\right),
    \end{split}\\
    \begin{split}
    \eta(4m_1\!+\!1,4m_2\!+\!2,4m_3\!+\!2) ={}& \myfrac{1}{4}\big((-1)^{m_2+m_3}\eta(m_1,m_2,m_3)\!+\! (-1)^M\theta(m_1,m_2,m_3)  \\ 
    & \hspace{3cm} \left. -\!(-1)^{m_1}\eta(m_1,m_2\!+\!1,m_3\!+\!1) \!-\!\theta(m_1\!+\!1,m_2\!+\!1,m_3\!+\!1)\right),
    \end{split}\\
    \begin{split}
    \eta(4m_1\!+\!1,4m_2\!+\!2,4m_3\!+\!3) ={}& \myfrac{1}{4}\left(-\!(-1)^{m_2+m_3}\eta(m_1,m_2,m_3)\!-\! (-1)^{m_1+m_2}\eta(m_1,m_2,m_3\!+\!1) \right. \\ 
    & \hspace{3cm} \left. -\!(-1)^{m_1}\eta(m_1,m_2\!+\!1,m_3\!+\!1) \!+\!(-1)^{m_3}\eta(m_1\!+\!1,m_2\!+\!1,m_3\!+\!1)\right),
    \end{split}\\
    \begin{split}
    \eta(4m_1\!+\!1,4m_2\!+\!3,4m_3\!+\!3) ={}& \myfrac{1}{4}\left((-1)^{m_2+m_3}\eta(m_1,m_2,m_3)\!-\! (-1)^{m_1}\eta(m_1,m_2\!+\!1,m_3\!+\!1) \right. \\ 
    & \hspace{3cm} \left. +\!(-1)^{m_1}\theta(m_1,m_2\!+\!1,m_3\!+\!1) \!-\!(-1)^{m_2+m_3}\theta(m_1\!+\!1,m_2\!+\!1,m_3\!+\!1)\right),
    \end{split}\\
    \eta(4m_1\!+\!2,4m_2\!+\!2,4m_3\!+\!2) ={}& 0,\\
    \begin{split}
    \eta(4m_1\!+\!2,4m_2\!+\!2,4m_3\!+\!3) ={}& \myfrac{1}{4}\big(-(-1)^M\theta(m_1,m_2,m_3)\!+\! (-1)^{m_1+m_2}\eta(m_1,m_2,m_3\!+\!1)  \\ 
    & \hspace{3cm} \left. +\!(-1)^{m_3}\eta(m_1\!+\!1,m_2\!+\!1,m_3\!+\!1) \!+\!\theta(m_1\!+\!1,m_2\!+\!1,m_3\!+\!1)\right),
    \end{split}\\
    \begin{split}
    \eta(4m_1\!+\!2,4m_2\!+\!3,4m_3\!+\!3) ={}& \myfrac{1}{4}\big((-1)^M \theta(m_1,m_2,m_3)  
        \!-\!(-1)^{m_1}\theta(m_1,m_2\!+\!1,m_3\!+\!1) \\
        & \hspace{3cm} \left. +\!(1\!-\!(-1)^{m_2+m_3})\theta(m_1\!+\!1,m_2\!+\!1,m_3\!+\!1)\right),
    \end{split}\\
    \eta(4m_1\!+\!3,4m_2\!+\!3,4m_3\!+\!3) ={}& \myfrac{1}{4}\bigl(\!-\!(-1)^M\theta(m_1,m_2,m_3)  
        \!+\!(1\!+\!(-1)^M)\eta(m_1\!+\!1,m_2\!+\!1,m_3\!+\!1)\!+\theta(m_1\!+\!1,m_2\!+\!1,m_3\!+\!1)\bigr),\\
\end{align*}
}

{\footnotesize
\begin{align*}
    \theta(4m_1,4m_2,4m_3) ={}& 0,\\
    \theta(4m_1,4m_2,4m_3\!+\!1) ={}& \myfrac{1}{4}\bigl((1\!-\!(-1)^M) \eta(m_1,m_2,m_3)  
        \!-\!(-1)^{m_3}\theta(m_1,m_2,m_3)\!+\!(-1)^{m_1+m_2}\theta(m_1,m_2,m_3\!+\!1)\bigr),\\
    \theta(4m_1,4m_2,4m_3\!+\!2) ={}& \myfrac{1}{2} \bigl((-1)^{m_3}\theta(m_1,m_2,m_3)+(-1)^{m_1+m_2}\theta(m_1,m_2,m_3\!+\!1) \bigr),\\
    \theta(4m_1,4m_2,4m_3\!+\!3) ={}& \myfrac{1}{4}\bigl(-\!(-1)^{m_3} \theta(m_1,m_2,m_3)\!-\!(1\!+\!(-1)^M)
        \eta(m_1,m_2,m_3\!+\!1) \!+\!(-1)^{m_1+m_2}\theta(m_1,m_2,m_3\!+\!1)\bigr),\\
    \theta(4m_1,4m_2\!+\!1,4m_3\!+\!1) ={}& \myfrac{1}{4}\bigl((1\!-\!(-1)^{m_2+m_3}\!+\!(-1)^M)
        \eta(m_1,m_2,m_3) \!-\!(-1)^{m_1}\eta(m_1,m_2\!+\!1,m_3\!+\!1)\bigr),\\
    \begin{split}
    \theta(4m_1,4m_2\!+\!1,4m_3\!+\!2) ={}& \myfrac{1}{4}\left(\!(-1)^{m_2+m_3}\eta(m_1,m_2,m_3)\!+\! (-1)^{m_3}\theta(m_1,m_2,m_3) \right. \\ 
    & \hspace{3cm} \left. -\!(-1)^{m_1+m_2}\theta(m_1,m_2,m_3\!+\!1) \!-\!(-1)^{m_1}\eta(m_1,m_2\!+\!1,m_3\!+\!1)\right),
    \end{split}\\
    \begin{split}
    \theta(4m_1,4m_2\!+\!1,4m_3\!+\!3) ={}& \myfrac{1}{4}\left( -\!(-1)^{m_3}\theta(m_1,m_2,m_3)\!-\!(-1)^{m_2}(1\!+\!(-1)^{m_1})\theta(m_1,m_2,m_3\!+\!1) \right. \\
    & \hspace{3cm} \left. +\!(-1)^{m_1+m_3}\theta(m_1,m_2\!+\!1,m_3\!+\!1)\right),
    \end{split}\\
    \theta(4m_1,4m_2\!+\!2,4m_3\!+\!2) ={}& 0,\\
    \begin{split}
    \theta(4m_1,4m_2\!+\!2,4m_3\!+\!3) ={}& \myfrac{1}{4}\left(-\!(-1)^{m_2+m_3}\eta(m_1,m_2,m_3)\!+\! (-1)^{m_2}\theta(m_1,m_2,m_3\!+\!1) \right. \\ 
    & \hspace{3cm} \left. +\!(-1)^{m_1}\eta(m_1,m_2\!+\!1,m_3\!+\!1) \!+\!(-1)^{m_1+m_3}\theta(m_1,m_2\!+\!1,m_3\!+\!1)\right),
    \end{split}\\
    \theta(4m_1,4m_2\!+\!3,4m_3\!+\!3) ={}& \myfrac{1}{4}\bigl((-1)^{m_2+m_3}
        \eta(m_1,m_2,m_3) \!+\!(-1\!+\!(-1)^{m_1}\!-\!(-1)^M)\eta(m_1,m_2\!+\!1,m_3\!+\!1)\bigr),\\
    \theta(4m_1\!+\!1,4m_2\!+\!1,4m_3\!+\!1) ={}& \myfrac{1}{4}\bigl((1\!-\!(-1)^M) \eta(m_1,m_2,m_3)  
        \!-\!(-1)^{M}\theta(m_1,m_2,m_3)\!+\theta(m_1\!+\!1,m_2\!+\!1,m_3\!+\!1)\bigr),\\
    \begin{split}
    \theta(4m_1\!+\!1,4m_2\!+\!1,4m_3\!+\!2) ={}& \myfrac{1}{4}\big( ((-1)^{m_3}\!+\!(-1)^M)\theta(m_1,m_2,m_3)\!+\!(-1)^{m_1+m_2}\theta(m_1,m_2,m_3\!+\!1)  \\
    & \hspace{3cm} \left. +\! \theta(m_1\!+\!1,m_2\!+\!1,m_3\!+\!1)\right),
    \end{split}\\
    \begin{split}
    \theta(4m_1\!+\!1,4m_2\!+\!1,4m_3\!+\!3) ={}& \myfrac{1}{4}\big(-\!(-1)^{m_3}\theta(m_1,m_2,m_3)\!-\! (-1)^{m_1+m_2}\eta(m_1,m_2,m_3\!+\!1)  \\ 
    & \hspace{3cm} \left. +\!(-1)^{m_1+m_2}\theta(m_1,m_2,m_3\!+\!1) \!-\!(-1)^{m_3}\eta(m_1\!+\!1,m_2\!+\!1,m_3\!+\!1)\right),
    \end{split}\\
    \begin{split}
    \theta(4m_1\!+\!1,4m_2\!+\!2,4m_3\!+\!2) ={}& \myfrac{1}{4}\big((-1)^{m_2+m_3}\eta(m_1,m_2,m_3)\!-\! (-1)^M\theta(m_1,m_2,m_3)  \\ 
    & \hspace{3cm} \left. -\!(-1)^{m_1}\eta(m_1,m_2\!+\!1,m_3\!+\!1) \!+\!\theta(m_1\!+\!1,m_2\!+\!1,m_3\!+\!1)\right),
    \end{split}\\
    \begin{split}
    \theta(4m_1\!+\!1,4m_2\!+\!2,4m_3\!+\!3) ={}&  \myfrac{1}{4}\left(-\!(-1)^{m_2+m_3}\eta(m_1,m_2,m_3)\!+\! (-1)^{m_1+m_2}\eta(m_1,m_2,m_3\!+\!1) \right. \\ 
    & \hspace{3cm} \left. -\!(-1)^{m_1}\eta(m_1,m_2\!+\!1,m_3\!+\!1) \!-\!(-1)^{m_3}\eta(m_1\!+\!1,m_2\!+\!1,m_3\!+\!1)\right),
    \end{split}\\
    \begin{split}
    \theta(4m_1\!+\!1,4m_2\!+\!3,4m_3\!+\!3) ={}& \myfrac{1}{4}\left((-1)^{m_2+m_3}\eta(m_1,m_2,m_3)\!-\! (-1)^{m_1}\eta(m_1,m_2\!+\!1,m_3\!+\!1) \right. \\ 
    & \hspace{3cm} \left. -\!(-1)^{m_1}\theta(m_1,m_2\!+\!1,m_3\!+\!1) \!+\!(-1)^{m_2+m_3}\theta(m_1\!+\!1,m_2\!+\!1,m_3\!+\!1)\right),
    \end{split}\\
    \theta(4m_1\!+\!2,4m_2\!+\!2,4m_3\!+\!2) ={}& \myfrac{1}{2} \bigl((-1)^M\theta(m_1,m_2,m_3)+\theta(m_1\!+\!1,m_2\!+\!1,m_3\!+\!1) \bigr),\\
    \begin{split}
    \theta(4m_1\!+\!2,4m_2\!+\!2,4m_3\!+\!3) ={}& \myfrac{1}{4}\big(-(-1)^M\theta(m_1,m_2,m_3)\!-\! (-1)^{m_1+m_2}\eta(m_1,m_2,m_3\!+\!1)  \\ 
    & \hspace{3cm} \left. -\!(-1)^{m_3}\eta(m_1\!+\!1,m_2\!+\!1,m_3\!+\!1) \!+\!\theta(m_1\!+\!1,m_2\!+\!1,m_3\!+\!1)\right),
    \end{split}\\
    \begin{split}
    \theta(4m_1\!+\!2,4m_2\!+\!3,4m_3\!+\!3) ={}& \myfrac{1}{4}\big((-1)^M \theta(m_1,m_2,m_3)  
        \!+\!(-1)^{m_1}\theta(m_1,m_2\!+\!1,m_3\!+\!1) \\
        & \hspace{3cm} \left. +\!(1\!-\!(-1)^{m_2+m_3})\theta(m_1\!+\!1,m_2\!+\!1,m_3\!+\!1)\right),
    \end{split}\\
    \theta(4m_1\!+\!3,4m_2\!+\!3,4m_3\!+\!3) ={}& \myfrac{1}{4}\bigl(\!-\!(-1)^M\theta(m_1,m_2,m_3)  
        \!-\!(1\!+\!(-1)^M)\eta(m_1\!+\!1,m_2\!+\!1,m_3\!+\!1)\!+\theta(m_1\!+\!1,m_2\!+\!1,m_3\!+\!1)\bigr).\\
\end{align*}
}

\end{document}